\documentclass[12pt]{article}
\usepackage{amscd}
\usepackage{amssymb}
\usepackage{amsthm}
\usepackage{amsmath}
\usepackage{graphics}

\newcommand{\zz}{ {\zeta} }
\newcommand{\bigo}{ {\mathcal{O}} }
\newcommand{\littleo}{ {\mathrm{o}} }
\newcommand{\binomial}[2]{\left(\begin{array}{c}{#1}\\{#2}\end{array}\right)}

\newcommand{\vol}{\mathbf Vol\ }

\newcommand{\bydef}{\stackrel{\mbox{\scriptsize def}}{=}}

\newtheorem {theorem}     {Theorem}
\newtheorem*{theoquo}     {Theorem}
\newtheorem {lemma}     {Lemma}
\newtheorem*{lemma*}     {Lemma}
\newtheorem {proposition}     {Proposition}
\theoremstyle{definition}
\newtheorem {definition}{Definition}
\newtheorem {axiom}{Axiom}

\newtheorem {remark}{Remark}
\newtheorem {example}{Example}

\newlength{\gnat}
\newcommand{\newsquare}[1]{{\square \hspace{-\gnat} 
       {\makebox[\gnat]{\raisebox{.25ex}{#1}}}}}

\newcommand{\renplus}[1]{{\boxplus_{#1}}}

\newcommand{\rentimes}[1]{\boxtimes}
\newcommand{\rendiv}[1]{{\newsquare{:}}}
\newcommand{\renpow}[2]{\ ^{\fbox{\tiny $#2$}}}
\newcommand{\renscal}[2]{\fbox{\tiny $#2$}_{#1}}
\newcommand{\bigrenplus}[1]{{\fbox{+}}_{#1}}
\newcommand{\RT}[2]{{\mathbb R}^{#1} \times {\mathbb T}^{#2} }
\newcommand{\orb}{{\mathrm{orb}}}
\newcommand{\smthg}{{ Lipschitz with constant $C_{k}$ }}
\newcommand{\dist}{{ \mathrm{d}}}

\newcommand{\newtext}[2]{{#1}}
\newcommand{\oldtext}[1]{}

\title{On the Geometry of Graeffe Iteration}

\author{Gregorio Malajovich\\
        {\small \em Dep. de Matem\'atica Aplicada, 
	Universidade Federal do Rio de Janeiro}\\
	{\small \em Caixa Postal 68530,
	Rio de Janeiro, RJ, 21945 -- BRASIL}\\
	{\small \em gregorio@labma.ufrj.br, http://www.labma.ufrj.br/\~{}gregorio}
        \\ 
        \and 
        Jorge P. Zubelli\\
        {\small \em IMPA}\\ 
	{\small \em Estrada Dona Castorina 110,
	Rio de Janeiro, RJ, 22460-320, BRASIL}\\
	{\small \em zubelli@impa.br, http://www.impa.br/\~{}zubelli}\\}

\date{Second revision, Aug 26, 1999}

\begin{document}
\maketitle
    \begin{abstract}
	A new version of the Graeffe algorithm for finding all
	the roots of univariate complex polynomials is 
	proposed. 
	It is obtained from the classical algorithm by a process
	analogous to renormalization of dynamical systems.
        \par
        This iteration is called Renormalized Graeffe Iteration.
        It is globally convergent,
        with probability 1. All quantities involved in the 
        computation are bounded, once the initial polynomial is
        given (with probability 1). This implies 
        remarkable stability properties for the new algorithm,
        thus overcoming known limitations of the classical
        Graeffe algorithm.
	\par
	If we start with a degree-$d$ polynomial, each renormalized
	Graeffe iteration costs $\bigo(d^2)$ arithmetic operations, with
	memory $\bigo(d)$. 
	\par
	A probabilistic global complexity bound is given. The case
	of univariate real polynomials is briefly discussed.
        \par
        A numerical implementation of the algorithm presented herein
	allowed us to solve random polynomials of degree up to 1000.
    \end{abstract}

\newpage
\vspace{\stretch{1}}
\ \\

\tableofcontents

\vspace{\stretch{1}}
\newpage
\section{Introduction}
\label{Sec1}
	\par	{\em Graeffe Iteration} was \oldtext{certainly}
		\newtext{one of}{A6\S 2(2)} the most prestigious 
		XIX century algorithm for finding roots of polynomials.
		At that time, computations were performed by hand, by people payed 
		specifically to perform those computations. They
		were called {\em calculateurs}~\cite{OSTROWSKII} or
                {\em computers}~\cite{BABBAGE}.
	\par    Let $f$ be a univariate polynomial, of degree $d$. Its
		Graeffe iterate is defined as:
	\begin{equation}
		Gf(x) = (-1)^d f(\sqrt{x}) f(-\sqrt{x}) 
	\end{equation}	
	\par 	This defines a many-to-one mapping in the space
		of all degree-$d$ polynomials (real or complex,
		as wish). The effect of this mapping is to
		square each root of $f$.
	\par 	After a few Graeffe iterations, the roots of $G^k f$ have
		(hopefully) incommensurate moduli. This is not true for
		\newtext{complex-}{A5(1)}conjugate roots, 
		which can be worked out in a different way.
	\par 	Assume, for simplicity, that $f$ is a complex polynomial,
		with no two roots of the same modulus. Then, $G^k f$ can
		be written
	\begin{equation}
		G^k f(x) = \sum_{i=0}^d a^{(k)}_i x^i
	\end{equation}
		where $a^{(k)}_i$ is given by the \oldtext{$d-i$-symmetric}
		\newtext{$(d-i)$-th symmetric}{A7(1)} function
		of the roots of $G^k f$. Therefore, $a^{(k)}_i$ is dominated
		by:
	\begin{equation}
		(-1)^{d-i} {\zeta_1}^{2^k} {\zeta_2}^{2^k} 
		\dots {\zeta_{d-i}}^{2^k}
	\end{equation}
		where $\zeta_1$, \dots, $\zeta_d$ are the roots of $f$
		ordered with decreasing modulus. (A more rigorous 
		statement of this will appear in Section~
		\ref{Sec4})
		\par
		Therefore, $-a^{(k)}_{i-1} / a^{(k)}_{i}$ is a good approximation
		for ${\zeta_{d-i+1}}^{2^k}$.
		\par
		\oldtext{Hence, we can easily obtain $|\zeta_i|$ for all $i$.
		There are many classical algorithms to recover the
		actual value of $\zeta_i$.}
		\newtext{Hence it is computationally easy to approximate
		$|\zeta_i|$ for all $i$. Although we also obtain
		$\arg \zeta_i \mod 2^{1-k} \pi$, we will discard this
		information in this paper to avoid additional complications.
		There are many classical algorithms
		to recover the actual value of $\zeta_i$. 
		}{A5(2-3)}
		See Pan~\cite{PAN} for
		a discussion. 
		\footnote{Added in the revised version: We deal
		with this issue in ~\cite{TANGRA}.}
	\medskip
	\par 	
	\oldtext{In this note, we claim that Graeffe iteration
		(after a suitable change of coordinates) is
		a globally convergent algorithm with probability $1$
		(Theorem~1 below).}
	\newtext{In this note, we apply a suitable non-uniform
	        change of coordinates (renormalization)
		to the Graeffe iteration 
		operator, to make it `convergent' with probability $1$.
\par
		The algorithm obtained by this change of coordinates will
		be called {\em Renormalized Graeffe Iteration}. 
		We use the following systems of coordinates for
		each iterate $G^{k} f$ of the Graeffe method applied to
		$f$:
		% $f$, $Gf$, $G^2f$, $\cdots$, $G^k(f)$, $\cdots$:
		The coefficients of $G^k(f)$ and all related intermediate
		computations will be represented in scaled polar coordinates,
		where a complex number $w$ is represented by `magnitude'
		$2^{-k} \log_2 |w|$ and `argument' $\arg w \in [-\pi, \pi]$.
		Calculations will always be performed `in coordinates'. 
		The `magnitude' variables of $G^k f$ `in coordinates' will
		converge with probability~1 (Theorem~\ref{th1} below).
\par
		The precise construction of the Renormalized Graeffe Operator
		is postponed to Section~\ref{Sec2}.
	\medskip
}{A2\S -1, A5(4-5) and B\S 3}
	\par    We also claim that the Renormalized Graeffe algorithm 
		compares well with
		available numerical software or theoretical algorithms.
\newtext{\par
        However, in this paper, we are considering a modified problem:
	our algorithm is designed to find the absolute values of the
	roots, not the actual roots. Therefore we will compare its
	complexity to the complexity of finding the absolute value
	of the roots by other existing algorithms. 
	\par
	In \cite{TANGRA}, we explain how to
	modify this algorithm to obtain the actual roots, without
	endangering the complexity estimates.
	\par
	} {A2(1)}
		The algorithm presented here has arithmetic complexity
		 $\bigo(d^{2})$ for each iteration, 
		and  memory size $\bigo(d)$, where $d$ is the degree
		of the polynomial.
%\newtext{The number of steps for obtaining an accurate value of the
%         roots is $O(\log d)$, see Theorem~\ref{th2} below.}{A2(2)}
	\par
		\oldtext{Complexity}\newtext{The number of iterations}{A2(2)}
                 will be bounded also in terms of a probability of
		failure (Theorem~\ref{th2} below). 
		\oldtext{This will be possible due to}
		\newtext{The authors personally believe that this is only
		possible due to}{A5(6)}
		the clean mathematical structure of Graeffe iteration.
		Our bound improves previous probabilistic bounds 
		\newtext{(in the sense of probability of success)}{A7(2)}
		on the
		complexity of solving polynomials. 
		(See Renegar~ \cite{RENEGAR} and Shub-Smale~\cite{BEZIV}).
	\par
		Also, our algorithm compares well with practical software,
		like for instance the algorithm
		in Matlab (running time $\bigo(d^{3})$ and memory $\bigo(d^{2})$).
		\newtext{Instead, Renormalized Graeffe Iteration
                seems to run in time $\bigo(d^2)$. It also}
		{A2(2)}
		\oldtext{Renormalized Graeffe iteration} seems much more stable,
		for \newtext{most (in a probabilistic sense)}{A3\S 1}
		\oldtext{generic}
		complex polynomials of degree, say, 1000. (See
		Section~\ref{Sec6} for a discussion of preliminary
                experimental results).
	\par
\newtext{In order to compare with deterministic algorithms one should
         bear in mind that the algorithm presented here is probabilistic
	 and may be quite slow on a set of non-zero measure. Also, 
	 the complexity of deterministic algorithms such as
	 \cite{NEFF-REIF,PAN,SCHONHAGE0,SCHONHAGE}
	 can be given in terms of number of arithmetic operations or
	 in number of bit operations. % The authors of this paper 
	 We believe 
	 that estimates of the order of  $\bigo(d^{1+\alpha})$, $\alpha >0$
	 in the number of operations are made at the expense of a
	 large increase in the working precision, of the order of
	 $\bigo(d^{2+\alpha})$ bits. However, this is not a lower bound
	 and it may as well happen that those algorithms can work
	 with substantially smaller precision on a large set of
	 inputs. 
\par
	In comparison, our Renormalized Graeffe Iteration was designed
	for floating-point arithmetic (our results assume an arbitrary,
	but fixed floating-point precision).}{A2(3--5)}
\oldtext{
		In order to compare it with theoretical algorithms
		of arithmetic complexity 
		$\bigo(d^{1+\alpha})$, $\alpha >0$, 
		(see \cite{NEFF-REIF,PAN,SCHONHAGE0,SCHONHAGE}) 
		one has to keep in mind that the bit complexity of
		those fast algorithms is (crudely) $\bigo(d^{3+\alpha})$, 
		since they 
		require at least $\bigo(d^{2+\alpha})$ bits of precision 
		(We are omitting here several factors). In comparison, our
		Renormalized Graeffe Iteration is, by construction, suitable 
		to usual floating-point arithmetic. 
}
		However, the details of the rounding-off
		analysis of renormalized Graeffe iteration will not be 
		discussed in this paper. 
	\par 	Instead, some preliminary 
		numerical results are presented in Section~\ref{Sec6}.
		They support the \newtext{empirical}{A5(7)}
		fact that \oldtext{generic}
		\newtext{typical random}{A3\S 1} degree 1000 polynomials
		can be solved within a precision of \oldtext{64}
		\newtext{64}{A5(7)} bits of mantissa \newtext{(IEEE 854 
		long double)}{A5(7)}. We
		expect a factor of a \oldtext{polylog of $d$}
		\newtext{polynomial in $\log d$}{A5(8)} 
		bits of mantissa to be
		\oldtext{necessary} \newtext{sufficient}{A5(8)} 
		in general, with probability 1. (This is a conjecture).
\newtext{
        \par    The result in Theorem~\ref{theo1} holds for ``reasonable''
                probability distributions on the space of polynomials. By
                reasonable, we mean all probabilities with
		bounded Radon-Nikodym derivative with respect to
                Lebesgue probability in the projectivization of
		coefficient space.
}{Changed.}

\begin{theorem} \label{th1}
                There is \oldtext{(we construct)}
		a renormalization of the Graeffe
		iteration, such that if $f$ is a \oldtext{generic}
		degree-$d$ polynomial
		(in a measure theoretical sense) then 
		\newtext{with probability 1}{A3\S 1}
		this renormalized Graeffe
		iteration produces $d+1$ sequences, each one converging to
		some $h_i$, s.t. $\log |\zeta_{i}| = h_i - h_{i+1}$, 
                and \oldtext{$\zeta_{i}$ is a root}
	        \newtext{$\zeta_1$, $\cdots$, $\zeta_d$ are roots}{A5(10)}	
		of $f$. 
		Moreover, each iteration can be performed in 
		$\bigo(d^2)$ arithmetical operations and all iterations can be 
		performed with memory $\bigo(d)$.
\label{theo1}
\end{theorem}
\par
\newtext{ This theorem is constructive, in the sense that an explicit 
          construction of the renormalized Graeffe
          iteration will be given.}{A5(9)}
	\par	In Section~\ref{Sec2}, we discuss the precise meaning of 
		renormalization in our context. 
		Its main consequence, will be to produce an
		algorithm operating on a bounded set of numbers. This solves
		the main stability problem of classical Graeffe iteration that
		prevented it from finding all roots at once. See for example
		Henrici's comments on \oldtext{Graeffe~\cite{HENRICI}}
		\newtext{FFT-based Graeffe iteration~\cite{HENRICI}
		(Vol~III, last paragraph of p. 69)}{A5(11)}. 
\oldtext{
\par    The definition of renormalization will outlaw, for instance,
		the following fast but unstable algorithm: Perform
		$k$ steps of Graeffe using FFT-based polynomial multiplication,
                where $k$ is given in Theorem~\ref{th2}
		}
\newtext{
   The definition of renormalization will outlaw FFT-based
   Graeffe iteration. Indeed, although FFT is known to be
   stable with respect to vector norms, it is not component-wise
   stable with respect to the relative error. This means that
   some of the coefficients of the $k$-th Graeffe iterate may have
   a large relative error, and hence some of the roots will be
   extremely inaccurate. This may be disastrous if one wants to 
   retrieve all the roots at the same time.}{B\S 2(1)}
   
\medskip
\par

% \par	The result in Theorem~\ref{theo1} holds for ``reasonable'' 
% probability distributions on the space of polynomials. By
% reasonable, we mean all absolutely continuous probability 
% distributions
% with respect to Lebesgue probability in coefficient space.
\newtext{
\begin{theorem}\label{th2}
Let $f$ be a random complex polynomial of degree $d$. Let $b \ge 1 + \log_2 d$.
Then, with probability $1-\delta$, $k$ steps of the Renormalized
Graeffe Iteration will approximate the $\log |\zeta_i|$'s with relative
precision $2^{-b}$, where
\[
k \ge c_1 + c_2 \log_2 b - c_3 \log_2 \delta
\]
and $c_2$, $c_3$ are universal constants. The constant $c_1$
depends on the choice of the probability distribution, and on
$d$.
\end{theorem}

Whenever speaking of random polynomials, we like to 
consider the normally invariant probability density introduced by
Kostlan~\cite{KOSTLAN} (See also Section~\ref{Sec5}). 
However, the above mentioned result is true for any reasonable 
probability distribution.

The experimental results in figure~\ref{results4} support the conjecture that,
under Kostlan's probability distribution, we can fix
\[
c_1 = c_4 \log_2 d
\]
where $c_4 \approx 2$ is a universal constant.

}{}
\par
\newtext{}{LARGE SECTION DELETED}
\medskip
\par
   We will briefly discuss the real case, and how to deal with
   \newtext{complex-}{A5(1)}conjugate 
   roots or roots with same modulus in Section ~\ref{Sec4}.

% from now on old

% and of bounded number
		% range will be explained in Section~\ref{explains}.
		% The model of computation is {\em idealized} IEEE
		% floating point arithmetic. This means we assume
		% that we have a large enough number of bits of mantissa 
		% and exponent. However, we do require this number of bits
		% to be finite, given any input. 
		% This model allows approximate and transcendental operations.  
	\medskip
	\par

	{\bf Historical remarks}. The Graeffe iteration was 
		developed independently by Dandelin (1826), by
		Graeffe (1837) and by Lobachevsky (1834). We call
		it Graeffe iteration to conform with most of
		the literature. See Householder~\cite{HOUSEHOLDER}
		for early references and priority questions.
		See Dedieu~\cite{DEDIEU} for an application of Graeffe
		\newtext{'s}{A7(4)}
		algorithm. 
% SEM A FOTO DO GRAEFFE E DO DANDELIN FICAMOS TENDENCIOSOS!
% \marginpar{ \special{psfile=lobatchevskii.ps voffset=-200 hoffset=-160 vscale=60 hscale=60 }
%	\ \\ {\tiny N.I.Lobachevsky}}
	\par 	Important theoretical results were obtained by
		Ostrowskii~\cite{OSTROWSKII} in 1940. Also, by that
		time, numerical analysis books mentioned Graeffe
		iteration as the preferred algorithm for 
                zero-finding (see e.g. Uspensky ~\cite{USPENSKII}
		\newtext{Page 318}{A6\S 2(1)}.
		For another early \oldtext{computational}
		\newtext{computer implemented}{A7(5)} algorithm, see 
		Bareiss~\cite{BAREISS1,BAREISS2} and also Blish and
		Curry~\cite{BC}).
	\par  	With the advent of digital computing, the practical
		use of Graeffe iteration seems to have been forgotten.

		Most popular zero-finding algorithms seem to be based
		now on QR iteration (Matlab) or in a several steps,
		root-finding
		plus deflation scheme. (e.g. Jenkins and Traub 
		\cite{JENKINS-TRAUB}). \oldtext{However} 
		\newtext{See, however,}{} Cardinal~\cite{CARDINAL},
		Edelman and Murakami~\cite{EDELMAN-MURAKAMI},
		Emiris, Galligo and Lombardi  ~\cite{EGH},
		and Toh and Trefethen~\cite{TREFETHEN-TOH}. 

	\par 	In a more theoretical perspective, Graeffe iteration is
		considered as a sort of pre-conditioning for polynomial
		splitting. Splitting a polynomial means factorizing it
		into one factor with \oldtext{very}\newtext{}{A6\S 2(1)}
		large roots, and another with
		\oldtext{very} small roots. 
		Splitting is used to obtain extremely
		fast theoretical algorithms (see Sch\"onhage~\cite{SCHONHAGE0,
		SCHONHAGE}, \newtext{Kirrinnis~\cite{KIRRINIS},
		Neff and Reif~\cite{NEFF-REIF},}{A7(6)}
		Bini and Pan~\cite{BP2},
		Mourrain and Pan~\cite{MP},
		Pan~\cite{PAN2,PAN3,PAN}, Pan et \oldtext{alli}
		\newtext{al.}{A7(7)}~\cite{PANETALLI}, 
		Malajovich and Zubelli
		~\cite{SPLITTING}). The main practical difficulty for those 
                algorithms seems to be the large precision required by 
		Graeffe iteration. Also, those algorithms are quite close
		to known lower bounds on \oldtext{arithmetic}
		\newtext{topological}{A5(15)} complexity 
                (See Vassiliev~\cite{VASSILIEV}). For a 
		related lower bound see Novak and Wo\'zniakowski~\cite{NW}. 
	\par 	An important paper by Grau in 1963~\cite{GRAU}
		laid some of the bases for a version of Graeffe iteration
		adapted to digital computers. He identified the problem
		of the increasing {\em numerical range}. During Graeffe
		iteration, some of the coefficients can become so large
		that the floating point system cannot accommodate them
		anymore.
	\par	While most of
		the literature suggests to find one root at a time and
		then use deflation, disregarding some stability problems
		(See e.g. Henrici~\cite{HENRICI}), Grau proposed a 
		globally convergent algorithm. Grau's algorithm would
		involve only bounded quantities. 
	\par	As far as we know, that
		paper was completely forgotten. The algorithm suggested
		by Grau has complexity $\bigo(d^2)$ and memory usage of
		$\bigo(d^2)$. It
		may be considered as the precursor of the one we 
		shall introduce below.

\section{Iterative Algorithms and Renormalization}
\label{Sec2}

	\par	In this paper, we will produce a version of Graeffe
		iteration that has {\em bounded} numerical range, for
		most input polynomials. The crucial concept in the
		construction of this algorithm is the idea of
		{\em renormalization.}
 	\medskip
	\par 	Renormalization is a tool used in understanding the qualitative
		behavior of iterative phenomena that range over
		different scales. A \oldtext{very}\newtext{}{A6\S 2(1)}
		rich theory of renormalization
		exists for one-dimensional dynamical systems.
		See Feigenbaum~\cite{FEIGENBAUM}, McMullen~\cite{MCMULLEN},
		and De Melo-Strien~\cite{MELO_STRIEN}. As for the 
		multi-dimensional case see Palis-Takens~\cite{G_D}.
	\par 	
	\oldtext{ Although we do not want to propose a general theory of
		renormalization of iterative algorithms, we may 
		illustrate what renormalization means by a well-known
		example and then give some tentative definitions.
		}
        \newtext{We shall illustrate what we mean by renormalization in
		the setting of iterative algorithms by a well-known
		example and then proceed with the definitions}{A5(16)}
	\par
\begin{example}
		Let $f: \mathbb C \rightarrow \mathbb C$ 
		be a smooth function. \oldtext{$\mathbb C \rightarrow \mathbb C$}
		Newton iteration associates to a point $x$ the point
		$N(x) = x - f(x) / f'(x)$. The sequence of Newton iterates 
		converges to a zero of $f$, provided that $x$ is picked
		close enough to a non-degenerate zero. 
	\par	We assume now that $f(x) = 0 + f'(0) x + \text{h.o.t.}$
		\newtext{}{B\S 4(3)},
		and that $f'(0) \ne 0$. We will consider now the 
                renormalization of the {\em iterative algorithm} (or
                mapping) $x \mapsto N(x)$. Although we are dealing here
	        with a simple example (the answer is always $0$),
	        its renormalization will have some of the main features
                of the renormalized Graeffe iteration, yet to be defined.
	\par 	\oldtext{Therefore, we should} 
		\newtext{The basic idea, therefore, is to}{for clarity} 
		look at Newton iteration with a variable
                {\em microscope} of variable lens. \oldtext{In our example we are
		performing Newton iteration $N$ of a univariate
		polynomial $f$ with a root located
		at the origin.}\newtext{}{A5(17)} We look at the mapping:
	\begin{equation}
		N_{\epsilon} = h_{\epsilon ^{-2}} 
				\circ N
				\circ h_{\epsilon }
	\end{equation} 
		where $h_{\epsilon}$ means the homothety of ratio $\epsilon$.
        \par 	When $\epsilon$ tends to zero, $N_{\epsilon}$ tends 
		to the map: $y \mapsto \gamma y^2$, defined from the disk 
		$D(|\gamma^{-1}|)$ of
		radius $|\gamma^{-1}|$ onto itself.
		Here, $\gamma$ is a parameter
		chosen equal to \oldtext{$f^{(2)}(0) / f'(0)$}
		\newtext{$f^{(2)}(0) / 2 f'(0)$}{B\S 4(4)}.	
		(Proof of this fact: Taylor's Theorem.
		The choice of the radius $|\gamma|^{-1}$ makes
		the limiting map surjective).
\end{example}
	\par 	This process (which we call renormalization) gives us
		qualitative information on the dynamics of Newton 
		iteration near a root. We can summarize this information
		in an {\em eventually commutative} diagram (commutative in the 
		limit):
\oldtext{
        \begin{equation}
	\begin{CD}
                D(|\gamma^{-1}|\epsilon) @>\left({h_{\epsilon}}\right)^{-1}
		>> D(|\gamma^{-1}|) \\
		@V{N}VV            @VV{y \mapsto \gamma y^2}V \\
                D(|\gamma^{-1}|\epsilon^{-2}) 
		@>{\left(h_{\epsilon} \right)}^2>> D(|\gamma^{-1}| )
	\end{CD}
	\end{equation}}
\newtext{
        \begin{equation}
	\begin{CD}
         D(|\gamma^{-1}|\epsilon) & \ni & x
	  @>\left({h_{\epsilon}}\right)^{-1} >> 
	  y & \in & D(|\gamma^{-1}|) \\
& &	@V{N}VV            @VV{y \mapsto \gamma y^2}V \\
      D(|\gamma^{-1}|\epsilon^{-2}) & \ni & N(x)
	@>{\left(h_{\epsilon} \right)}^{-2}>> 
		\gamma y^2 & \in & D(|\gamma^{-1}| )
	\end{CD}
	\end{equation}}{A3\S 3(1) and A6\S 1(18)}

	\par 	In the example above, the homothety $\left(h_{\epsilon}
		\right)^{-1}$ is
		generating the renormalization group. In general, 
		we want to consider renormalizations that lead to
		a commutative diagram, on the limit:
	\begin{equation}
	\begin{CD}
                y @>R>> R(y)\\
		@VGVV            @VV{G^R}V \\
                G(y) @>>R^2> R^2(G(y)) 
	\end{CD}
	\end{equation}
	\par 	Above, $G$ is the original algorithm (possibly after 
	\oldtext{eventually}\newtext{}{A7(8)}
		some uniform change of coordinates). We denote by
		$R$ the renormalization map, and by $G^R$ the renormalized
		version of our iteration. Usually, $G^R$ depends on
		a renormalization parameter. However, we want $G^R$ to 
		converge to a limiting map, with a simple dynamics. Moreover,
                computation of $G^R$ should be `stable', in a \oldtext{very}
		\newtext{}{A6\S 2(2)}
		precise sense. 
	\medskip
	\par 	
		This suggests the following heuristics, in the case
                of Graeffe iteration:
		we would like to consider
		a polynomial $f$ represented by its roots. More precisely,
		a polynomial $f$ can be represented by the vector
		$(\log |\zeta_i|, \arg(\zeta_i))$ $\in$ $\mathbb R^d \times
		\mathbb T^d$, where $\zeta_i$ is the $i$-th root of $f$,
		and $\mathbb T$ denotes the 
		additive group $\mathbb R / 2 \pi \mathbb Z$. In that case,
		Graeffe iteration is just multiplication by 2. 
	\par 	Renormalization would be a division by 2 of the
		log of the radii of the roots. We would have:
\oldtext{
\begin{equation}
%\label{CDlim}
	\begin{CD}
                \mathbb R^d \times \mathbb T^d 
		@>R^k>>
		\mathbb R^d \times \mathbb T^d  \\
	        @V{\times 2}VV 
		@VV(r,\theta)\mapsto(r,2\theta)V\\
                \mathbb R^d \times \mathbb T^d 
		@>>R^{k+1}> 
		\mathbb R^d \times \mathbb T^d 
	\end{CD}
	\end{equation}
}
\newtext{
\begin{equation}
\label{CDlim}
	\begin{CD}
                \mathbb R^d \times \mathbb T^d 
		& \ni & (\log |\zeta_i|, \arg \zeta_i)
		@>R^k>>
		(r,\theta) & \in & 
		\mathbb R^d \times \mathbb T^d  \\
	        & & @V{\times 2}VV 
		@VV(r,\theta)\mapsto(r,2\theta)V\\
                \mathbb R^d \times \mathbb T^d 
		& \ni & (2 \log |\zeta_i|, 2 \arg \zeta_i \mod 2\pi)
		@>>R^{k+1}> 
		(r, 2 \theta \mod 2 \pi) & \in &
		\mathbb R^d \times \mathbb T^d 
	\end{CD}
	\end{equation}
}{A3\S 3(1)}

\par 	Therefore, the limit of the renormalized map should be
                the map: $(r,\theta) \mapsto (r,2 \theta)$. 
                Of course, we do not know in advance the roots of the
		polynomial. Instead, we will produce a chain of
		commutative diagrams `converging' to diagram~(\ref{CDlim}).
		In order to do that,
	 	let us assume that a polynomial
	\begin{equation}
		a_0 + a_1 x + a_2 x^2 + \dots + a_d x^d
	\end{equation} 
		is represented by the vector
\newtext{
        \begin{equation}
		(\hat a_i, \alpha_i) =
		\left(
	        \log \left| a_i \right|,	
		\arg( a_i)
		\right)
        \end{equation} 
	}{A6\S 1(19)}
	\newtext{
	In this paper, for clarity of exposition, we will ignore
	the case where some of the $a_i$ is zero. We can do that
	because most of our results hold `with probability 1'. In
	practice, polynomials with a zero coefficient do arise.
	See Section ~\ref{Sec6} for implementation comments.
	}{B\S 4(5,6)}
	\par    
	        If the operator $G$ represents Graeffe iteration in this
		new system of coordinates, we will have:
\oldtext{
	\begin{equation}
% \label{cd}
        \begin{CD}
                @VVV                            
		@VVV \\
                G^k (f) \in \mathbb R^{d+1} \times \mathbb T^{d+1}  
		@>R^k>>
                (\hat a^k, \alpha^k) \in \mathbb R^{d+1} \times \mathbb T^{d+1} 
		\\
		@V{G}VV 
                @VV(G^{R,k})V        
		\\
		G^{k+1}(f) \in \mathbb R^{d+1} \times \mathbb T^{d+1} 
                @>>R^{k+1}> 
                (\hat a^{k+1}, \alpha^{k+1}) \in 
                      \mathbb R^{d+1} \times \mathbb T^{d+1} 
		\\
                @VVV                             
		@VVV
        \end{CD}
	\end{equation}
	}
\newtext{
	\begin{equation}
\label{cd}
        \begin{CD}
                & & @VVV                            
		@VVV \\
		\mathbb R^{d+1} \times \mathbb T^{d+1}  
                & \ni & G^k (f)  
		@>R^k>>
                (\hat a^k, \alpha^k) & \in & \mathbb R^{d+1} \times \mathbb T^{d+1} 
		\\
		& & @V{G}VV 
                @VV(G_k)V        
		\\
		\mathbb R^{d+1} \times \mathbb T^{d+1} 
		& \ni & G^{k+1}(f)  
                @>>R^{k+1}> 
                (\hat a^{k+1}, \alpha^{k+1}) & \in &
                      \mathbb R^{d+1} \times \mathbb T^{d+1} 
		\\
                & & @VVV                             
		@VVV
        \end{CD}
	\end{equation}
	}{A3\S 3(1,6)}
\newtext{
	\par
	Above, $\hat a^k = (\hat a_0^k, \cdots, \hat a_d^k)$ and
	$\alpha^k = (\alpha_0^k, \cdots,  \alpha_d^k)$. Also,
	$k$ is a superscript, not an exponent. However, 
	$R^{k}$ denotes the $k$-th iterate of $R$.}{A6\S 1(20)}
	\medskip
	\par	In this diagram, $R$ maps $(r,\theta)$ into $(r/2, \theta)$.
                Although this defines \oldtext{$G^{R,k}$}
		\newtext{$G_k$}{A3\S 3(6)} as a mapping, the 
		algorithmic construction of \oldtext{$G^{R,k}$} 
		\newtext{$G_k$}{A3\S 3(6)} is postponed to Section
		~\ref{Sec3}. On the limit, \oldtext{$G^{R,k}$} 
		\newtext{$G_k$}{A3\S 3(6)} `converges' to
		the mapping: \oldtext{$a, \theta \mapsto a, 2 \theta$}
		\newtext{$(a, \theta) \mapsto (a, 2 \theta)$}{B\S 4(7)}. 
		This will
		\oldtext{assure}\newtext{ensure}{A6\S 1(21)}
		that $\lim \hat a^k$ exists. We will 
		\oldtext{claim}\newtext{show}{A6\S 1(21)} that
		this limit satisfies \oldtext{(generically)}
		\newtext{with probability 1}{A3\S 1}
\[
		\lim_{k \rightarrow \infty} \hat a^k_j 
	      - \lim_{k \rightarrow \infty} \hat a^k_{j+1} 
		= \log |\zeta_{d-j}|
\]
 		where $\zeta_{d-j}$ is the \oldtext{$d-j$-th}
		\newtext{$(d-j)$-th}{A7(10)} root of the
		original polynomial, the roots being ordered by decreasing
		modulus..  
	\par 	Also, by using the renormalized Graeffe iteration
		\oldtext{$G^{R,k}$} 
		\newtext{$G_k$}{A3\S 3(6)} 
		instead of the classical Graeffe iteration
                $G$, we will 
		be able to bound all the intermediate calculations to
		a compact, depending on the input $f$. This will 
		imply several stability results.
	\medskip
	\par
		In order to define formally what we mean by a renormalized
	algorithm, we need to state the conditions that we expect
	\oldtext{$G^{R,k}$} 
	\newtext{$G_k$}{A3\S 3(6)} 
	to satisfy. These conditions will define a class
	of algorithms, that we call renormalized iterative algorithms.
	Some examples: usual Newton, renormalized Newton, all reasonable
	algorithms based on a contraction principle and, of course,
	Renormalized Graeffe Iteration (To be constructed). 
	\medskip
	\par

Before defining what we mean by a renormalized algorithm, we will
briefly discuss the notion of algorithm. There are many possible
definitions (See Blum, Cucker, Shub and Smale~\cite{BCSS}).

\begin{definition} \label{itealg}
An iterative algorithm $M$ is a Blum-Shub-Smale (BSS) machine
over ${\mathbb R}$, modified as follows: 
\begin{enumerate}
\item If the input is in ${\mathbb R}^{l}\times {\mathbb T}^{m}$  
 for a pair  \oldtext {$(l,k)$}
 \newtext{$(l,m)$}{A7(11) and B\S 4(8)} 
, then the output
is in ${\mathbb R}^{l}\times {\mathbb T}^{m}$. 
The integer $l+m$ is called the input size.
If the input is denoted by $x$, the output is denoted 
$M(x)$, and $M^k $ means the composition $M\circ M \circ \cdots 
$\newtext{$\circ$}{A7(12)}$M$
\newtext{\par Also, we say that the iterative algorithm `computes'
the function $x \mapsto \lim_{k \rightarrow \infty} M^k(x)$}{Avoid
confusion between iterate and limit.}
\item Computation nodes are allowed to perform 
\newtext{the following}{A4\S 1(1)} 
elementary functions\oldtext{, i.e.}\newtext{:}{}
polynomial evaluation, \newtext{absolute value,}{}
(real) logarithms, (real) exponential, sine, 
cosine, rational
functions, and compositions of those.

\item $M_{\epsilon}(x)$  will denote the result of approximating 
$M(x)$ by allowing each operation with non-integer parameters 
to be performed with relative precision $\epsilon$. (This is also
known in numerical analysis as the $(1+\epsilon)$ property.)
\newtext{The parameter $\epsilon$ is allowed to vary in $(0, \frac{1}{2})$.}
{A4\S 1(2)}
The machine is supposed to ``know" $\epsilon$. 
\newtext{This means that the value of $\epsilon$ can be used in intermediate
calculations.} {A4\S 1(2)}
Also, the approximation
is performed by some prescribed algorithm. (e.g., IEEE arithmetic
with $- \log_2 \epsilon$ bits of mantissa). 

\item \newtext{For all $\epsilon$,}{A4\S 1(2)}
the arithmetic complexity of $M$ applied to the input $x$ is the
number of \oldtext{arithmetic operations,}
elementary function evaluations 
\newtext{(from the list of item~1)}{A4\S 1(1)} and
branchings performed with input $x$. We require
the arithmetic complexity of the approximate algorithm 
\oldtext{
$M_{\epsilon}(x)$
to be the same as the arithmetic complexity of the original algorithm $M(x)$.
}
\newtext{
$M_{\epsilon}$ applied to the input $x$ 
to be the same as the arithmetic complexity of the original algorithm $M$
applied to input $x$.
}{A3\S 3(3)}
\end{enumerate}
\end{definition}

\begin{remark}
Item~1 will allow us to distinguish between real and angular variables,
the latter one being defined modulo $2\pi$. 
This will simplify notation when we speak of the distance between
points in ${\mathbb R}^{l}\times {\mathbb T}^{m}$. \newtext{Let $\dist$ be
that distance.}{A4\S 1(7)}
\end{remark}

\newtext{
\begin{remark}
	It has been argued that the outcome of a branching node would
	not be well-defined in the presence of numerical error. This is
	not true in the definition above. The branching nodes of 
	the machine $M$
	can be assumed, without loss of generality, 
	to branch on queries of the form $y>0$ or $y \ge 0$ or $y = 0$. 
	When the machine $M$ gets replaced by $M_{\epsilon}$, the
	branching nodes stay formally the same. The value of $y$, however, is
	contaminated with a certain numerical error. It is still
	a perfectly defined real number, and it may be compared to zero.
	Thus, the branching nodes branch `correctly', for a 
	slightly perturbed input.
	\par
	This would lead to disastrous results if an approximate machine
	enters a `loop' due to numerical errors. Item 4 in the
	definition above is there to preclude that sort of loop, by
	ensuring that the arithmetic complexity of each iteration
	does not depend on the working precision. \oldtext{ This is the case
	for all iterative algorithms that the authors know about.}
\end{remark}}{A4\S 1(3) and B\S 2(3,4)}
\newtext{
\begin{definition} A function $\varphi$ can be computed in {\em finite time}
if and only if there is a BSS machine over $\mathbb R$ 
modified as above items 2, 3 and 4,
that computes $\varphi$.
\end{definition}
}{Define the notion finite time precisley.}
One consequence of the previous definition is the following: suppose
a certain function $\varphi$ can be computed in finite time. Then, its branching
set is given by a finite set of equations. Outside the branching set,
$\varphi$ can be written (locally) as a composition of elementary functions.
Moreover, only a finite number of such compositions may appear,
one corresponding to each set of possible branchings.    

\medskip
\par
A few definitions \oldtext{that will be used in the sequel }
are in order now.  In the sequel we will need to use 
algorithms depending effectively on a parameter $k \in \mathbb N$.   
\oldtext{ Let $M$ be an iterative algorithm, with input $(f,k)$, 
where the parameter $k$ will keep track of iterations.
We write $M_{k}(f)=M(f,k)$.} \newtext{
        Those will
        be given by a machine $M$ with two inputs, say $k$ and $f$.
        However, we will denote the output as $M_k(f)$ and we
        will write $M_k$ for each of the input-output mappings
        obtained by restricting that machine to some fixed value
        of $k$.
        Also, in that situation we will speak explicitly of the
        algorithm $(M_k)_{k \in \mathbb N}$ or $M_k$ for short.
\par
}{A4\S 1(5)}

\par
The sequence $\left(  M_{k} \right)_{k=1}^{\infty}$ can be considered
\newtext{as a sequence of mappings}{Syntax.} 
from $\RT{l}{m}$ into itself.  
The orbit of $f$ by the sequence $\left(  M_{k} \right)_{k=1}^{\infty}$ 
is the set 
\[ \orb f \bydef \left( f, M_{1}(f),M_{2}\circ M_{1}(f),\dots \right) 
\subset \RT{l}{m} \mbox{ .} \]  

This set should be understood as the orbit of $f$ by the non-autonomous
dynamical system $M_k$ (Not the semi-group !)
The closure of $\orb f$ will be denoted by
$\overline{\orb f}$.

We shall say that a subset of $\RT{l}{m}$ \oldtext{is generic}
\newtext{has full-measure}{A3\S 1} if its complement is 
contained in a set of null measure. \newtext{We say that some property
is true almost everywhere (a.e.) if that property is true in a full-measure
set.}{}
We can now define a {\em renormalized algorithm} and make precise our
concept of renormalization:

\begin{definition} \label{defren2}
The algorithm \newtext{$(M_k)_{k \in \mathbb N}$ (or
$M_k$ for short)}{A4\S 1(5)}\oldtext{$M_k$}
is said to be a {\em renormalized iterative algorithm 
to compute} 
\oldtext{$\varphi(f): \RT{l}{m}  \rightarrow {\mathbb R}^{s}$}
\newtext{$\varphi : \RT{l}{m}  \rightarrow {\mathbb R}^{s}$, 
$f \mapsto \varphi(f)$}{A3\S 3(4)}, where $0\le s \le l$,
  if, and only if, it satisfies the Axioms 1 through 4 below. 
\end{definition}
\begin{axiom}[Consistency] 	% Axiom 1
	For almost every $f\in {\mathbb R}^{l} \times {\mathbb T}^m$ we have  
\oldtext{
\[ 
	\lim_{k\rightarrow \infty} 
	\pi \circ M_{k}\circ M_{k-1}\circ \cdots \circ M_{1}(f)
	= \varphi(f) 
\]
}
\newtext{
\[ 
	\lim_{k\rightarrow \infty} 
	\left( \pi \circ M_{k}\circ M_{k-1}\circ \cdots \circ M_{1} \right)(f)
	= \varphi(f) 
\]}{A3\S 3(5)}
	where $\pi$ is the projection of $\RT{l}{m}$ onto the first $s$ 
	coordinates of ${\mathbb R}^{l}$.
\end{axiom}

\begin{axiom}[Arithmetic Complexity] 	% Axiom 2
	The arithmetic complexity of \oldtext{$M_{k}(f)$}
	\newtext{$M_k$ with input $f$}{A3\S 3(4)} is bounded in terms of the
	size \newtext{$l+m$}{Clarity.} of $f$, independently of 
	$k$ and the coefficients of $f$.
\end{axiom}

\begin{axiom}[Propagation] 	% Axiom 3
	For \oldtext{generic $f\in \RT{l}{m}$ } 
	\newtext{almost every $f \in \RT{l}{m}$,}{A3\S 1}
	there exists a compact neighborhood
	$V \subset \RT{l}{m}$ 
	of $\orb f$ (under $\{ M_{k} \}$) 
	and \oldtext{$C \ge 1$}\newtext{C}{A4\S 1(6)}
	such that $M_{k}|_{V}$ is \smthg % with constant $C_{k}$
	and eventually $C_{k}<C$. 
\end{axiom}
\begin{axiom}[Stability] 	% Axiom 4
	For \oldtext{generic $f\in \RT{l}{m}$} 
	\newtext{almost every $f \in \RT{l}{m}$, }{A3\S 1}
	there exists a compact neighborhood
        $V \subset \RT{l}{m}$ of $\orb f$ (under $\{ M_{k} \}$)
	and $B\in {\mathbb R}$ such that $\forall g \in V$ and $\forall k$, 
\[
        \dist\left( 
M_{k,\epsilon} (g)
        \  , \   M_{k}(g) \right)  
        < \epsilon B \mbox{ .} 
\]
\end{axiom}

\medskip
\par
	Our concept of renormalized iterative algorithm subsumes
	several `reasonable' properties of iterative algorithms.
	Axiom 1 allows our algorithm to carry more information than
	what is actually required at output. Yet, we want our algorithm to
	produce a sequence converging to the expected result, for almost 
	every input.
	Axioms 3 and 4 will rule out unstable algorithms. 
	The idea behind Axiom 2 is that any honest iterative algorithm
	should have bounded arithmetic complexity for each iteration.
	This prevents the use of multiple precision arithmetic to
	obtain stability at the expense of (possibly exponentially
	many) extra arithmetic operations.

\medskip
\par 
	We shall now explore some consequences of the 
	definition of renormalized iterative algorithm  we proposed above.
	More specifically, our first 
	goal is to obtain  a (non-uniform) error bound on
	the result of iterating $k$ times a renormalized algorithm
	with precision $\epsilon$.

\begin{lemma} \label{lemE}
        Let $M$ be a renormalized iterative algorithm 
	\oldtext{and $f$ generic}.
        Then, \newtext{for almost every $f$ and}{A3\S 1} 
	for each $k$, \oldtext{ there is an $\epsilon>0$}
	we have that \newtext{for sufficiently small $\epsilon>0$}{} 
\[
        \dist  \left(
                 M_{k,\epsilon} \circ \cdots \circ M_{1,\epsilon} (f)
                ,   M_{k} \circ \cdots \circ M_{1} (f)
        \right) 
        < k A^k B \epsilon \mbox{ . } \]
	\newtext{Here, $A$ and $B$, 
	 depend on $f$, but not on $k$.}{} 
        \oldtext{Moreover, one can take}
	\newtext{In particular, that will be true for values of
	$\epsilon$ of the form}{} 
\[ \epsilon \le \frac{\rho}{kA^kB} \mbox{ ,} \]
	where $\rho$ is defined below.
\end{lemma}
As an immediate consequence of Lemma~\ref{lemE} we have that
\[
        \|  \left( \pi \circ 
                 M_{k,\epsilon} \circ \cdots \circ M_{1,\epsilon} )(f)
                - ( \pi \circ M_{k} \circ \cdots \circ M_{1} )(f)
        \right) \|_2
        < k A^k B \epsilon \mbox{ .}
\]

\begin{proof}[Proof of Lemma~\ref{lemE}]
We now describe the construction of the constants in the statement of
Lemma~\ref{lemE}.
Combining Axioms 3 and 4, there exists a compact neighborhood $W$ of $\orb f$
such that:
\begin{enumerate}
\item 
Every mapping $M_{k}$ is Lipschitz on $W$. Moreover, there
exists a constant $A$ such that for every $k$, the norm of the Lipschitz
constant  
of $M_{k}$ is uniformly bounded by 
\oldtext{$A = \sup C_k$}
\newtext{$A = \max \left( 1, \sup_k C_k \right)$,}{A4\S 1(6)} 
$C_k$ as in Axiom~3. 
\item 
\oldtext{For $g\in W$, there exists}
\newtext{There exists}{A3\S 2(1)}
$\epsilon_{0}$ such that
$\forall \epsilon < \epsilon_{0}$ we have
\[ {\dist}(M_{k,\epsilon}(g),M_{k}(g)) < \epsilon B
\mbox{ ,} \hspace{1cm} \forall g \in W \mbox{ .} \] 

We then define $\rho$ as the minimum of $\epsilon_{0}$ and
the distance of $\orb f$ to the boundary
of $W$.
\end{enumerate}

We may now conclude the proof of Lemma~\ref{lemE}.
Set 
\[ g=(M_{k-1}\circ \dots \circ M_{1})(f) \mbox{ ,} \]
and
\[ h=(M_{k-1,\epsilon} \circ \dots \circ M_{1,\epsilon})(f) \mbox{ .} \]
In that case, we want to bound
\[ \dist(M_{k,\epsilon}(h),M_{k}(g)) \le \dist (M_{k,\epsilon}(h),M_{k}(h))
+ \dist(M_{k}(h),M_{k}(g)) \]
If $h\in W$, one can bound 
\[ \dist(M_{k}(h),M_{k,\epsilon}(h)) \le B\epsilon \]
and
\[ \dist(M_{k}(h),M_{k}(g)) \le A \ \dist(g,h) \mbox{ .}  \]
By induction, 
\[ \dist(g,h)  \le (k-1) A^{k-1} B \epsilon \mbox{ ,} \]
and hence
\[ \dist(M_{k,\epsilon} (h), M_{k}(g) )\le k A^{k} B \epsilon \mbox{ .} \]
The condition 
\[ \epsilon < \frac{\rho}{k A^{k} B} \]
guarantees that $h \in W$. 
\end{proof}

\medskip
\par
	If we use a Turing machine model (or any other classical
	discrete complexity model), we can perform all the operations of
	$M$ in finite precision, and obtain:

\begin{lemma} \label{lemT}
Let $M$ be a renormalized iterative algorithm.
\oldtext{ Fix a generic} \newtext{For almost every $f$,
assume that the truncation error is bounded by}{A3\S 1}
\oldtext{
\[ \| \pi \left(
                 M_{k} \circ \cdots \circ M_{1}  \right) (f)
                - \varphi(f)
        \|_2  \]
	}
\newtext{
\[ \| \left( \pi \circ
                 M_{k} \circ \cdots \circ M_{1}  \right) (f)
                - \varphi(f)
        \|_2   < E^{2^{k}} \mbox{ ,} \] 
where $ E = E(f) \in (0,1) $.}{A3\S 3(5)}
\oldtext{ is bounded by $E^{2^{k}}$}  
Then, the \oldtext{Turing} complexity of approximating $\varphi(f)$ with
precision $\delta$ is $\bigo((\log_{2} \frac{1}{\delta})^{1+\alpha})$
where $\alpha > 0 $ is arbitrarily \newtext{
small.}{}
\end{lemma}
\newtext{
It is assumed above that the cost of arithmetic with $l$ bits
of mantissa is $O(l^{1+\alpha})$, for all $\alpha > 0$. 
See~\cite{BP} pp~78-79 for a sharper bound on the complexity
of long integer multiplication.}{A6\S 1(22)}

\begin{proof}[Proof of Lemma~\ref{lemT}]
\par
	Let 
\[
	k = \left\lceil \log_2 \frac{-\log_2 \frac{\delta}{2}}{-\log_2 E} 
	\right\rceil
	\mbox{ ,}
\]
	so that 
\[
	E^{2^k} \le \frac{\delta}{2} \mbox{ .}
\]

	Choose $\epsilon$ so that Lemma~\ref{lemE} holds, i.e.,
	$\epsilon < \rho/k A^k B$, and such that
	$k A^k B \epsilon < \delta/2$. This can be done
	for
\oldtext{
\[
	\epsilon \in
	\littleo 
	\left(
	    \frac{\delta}{\log_2 \log_2 \delta^{-1} (\log_2 \delta^{-1})^
		{\log_2 A}}
	\right)
        \subset	
	\littleo (\delta^{1 + \alpha})
\]
}
\newtext{
\[
\epsilon < c_1 \delta (\log \delta^{-1})^{c_2}
\]
for some constants $c_1$ and $c_2$ dependent of $f$.
}{}
	The total cost of computing $M_k$ is 
        $\bigo \left( a (\log_2 \epsilon^{-1})^{1 + \alpha} \right)$,
	where $a$ is the arithmetic complexity of the algorithm and
	$\alpha$ is arbitrarily small. \oldtext{Hence,}
	\newtext{Since $f$ is fixed, $a$ is a constant and 
	$k < \log_2 \log_2 \delta^{-1}$}{}
	 the total cost is
\[
	\bigo \left( (\log_2 \delta^{-1})^{1+\alpha}\right)
\]
\newtext{for all $\alpha$ arbitrarily small.}{}
\end{proof}

\medskip
\par	The complexity bound of Lemma~\ref{lemT} is non-uniform in $f$.
	It is also non-effective, in the sense that we give no procedure
	to estimate $k$ and \oldtext{$m$}\newtext{$\epsilon$}{}
	without the knowledge of $\rho$, 
	$B$ and $A$. Indeed, those quantities may depend on $f$.
	It \newtext{would be of some interest}{} to bound those 
	quantities in a probabilistic
	setting, similar to Theorem~\ref{th2}.

\begin{definition}
	Let $M$ be an iterative algorithm to compute \oldtext{$\Phi(F)$}
	\newtext{$\Phi: F \mapsto \Phi(F)$,}{A3\S 3(4)} 
\newtext{i.e., $\lim_{n \rightarrow \infty} M^n(F) = \Phi(F)$).
        }{A4\S1 (8)}	
	We say
	that $M_k$ is a renormalization of $M$ if 
\begin{enumerate}
	\item $M_k$ is a renormalized iterative algorithm to compute $\varphi(f)$.
	\item There are functions $\psi$ and $\eta$, defined almost everywhere
        and computable in finite time,
	such that the diagram
\[
\begin{CD}
	F @>\psi>> \psi(F) \\
@V{\Phi}VV	@VV{\varphi}V\\
	\Phi(F) @<\eta<< \varphi(\psi(F))
        \end{CD}
\]
	commutes.
	\item There is a function $R$, 
		computable in finite time, so that the diagrams:
\[
        \begin{CD}
                F @>R^k \circ \psi>> R^k(\psi(F)) \\
                @V{M}VV        @VV{M_k}V \\
                M(F) @>R^{k+1} \circ \psi>> R^{k+1}(\psi(M(F))) \\
        \end{CD}
\]
	commute.
\end{enumerate}
\end{definition} 	

	Theorem~\ref{th1} can be stated now in a more concise way:
	Let \oldtext{$z(f)$}\newtext{$\zz: f \mapsto \zz(f)$}{A3\S 3(4)}
	be the function that associates, to any univariate
	degree $d$ polynomial $f$, its roots $\zeta_1 , \dots , \zeta_d$ ordered by
	decreasing modulus. We have the following: 

\begin{theoquo}
                There is \oldtext{(we construct)}
		a renormalization of the Graeffe
		iteration to compute $|\zeta(f)|=\left( |\zeta_1(f)|, 
                \dots, |\zeta_d(f)|\right)$.
		Each iteration has arithmetic complexity 
		 $\bigo(d^2)$  and needs
		memory $\bigo(d)$.
\end{theoquo}

\section{Recurrence Relations and the Renormalization of Graeffe}
\label{Sec3}

		It is time to construct the renormalized Graeffe
		\oldtext{operator} 
		\newtext{iteration}{A7(3)}.
		Let $f(x) = f_0 + f_1 x + \dots + f_d x^d$. Let $h = Gf$ be
		its Graeffe iterate. The coefficients of $h$ can be written
		as~:
	\begin{equation}
\label{g1}
		h_i = (-1)^d \sum
			_{\substack{
				0\le i-j \le d \\
				0\le i+j \le d \\ }}
			 (-1)^{i-j} f_{i-j} f_{i+j}
	\end{equation} 
	\par	For convenience, we rewrite (\ref{g1}) as: 
	\begin{equation}
\label{g2}
		h_i = (-1)^{d+i} {f_{i}}^2 + 
			2 \sum
			_{ 1\le j \le \min(i,d-i)}
			 (-1)^{d+i-j} f_{i-j} f_{i+j}
	\end{equation}

	\medskip
	\par	The next step is to write those equations in terms
		of the log of the coefficients. More precisely, we will
		have to deal with the following two quantities:
	\begin{equation}
\label{g7}
		f^{\log}_i \bydef \log \left| f_i \right|
	\end{equation}  
	\begin{equation}
\label{g8}
		f^{\arg}_i \bydef \arg f_i 
	\end{equation}  
	\medskip
	\par	It is possible now to construct the renormalized Graeffe
		\oldtext{operator} 
		\newtext{iteration}{A7(3)}
		\oldtext{$G^{R,k}$} 
		\newtext{$G_k$}{A3\S 3(6)}. 
		Recall from diagram (\ref{cd}) that 
		this \oldtext{operator} 
		\newtext{iteration}{A7(3)} will map 
		$2^{-k} f^{\log}$ and $f^{\arg}$ into $2^{-k-1} h^{\log}$
		and $h^{\arg}$.
	\par	For that purpose, we introduce the notation:
	\par
%JPZ	\centerline{
%JPZ	\begin{tabular}{ll}
%JPZ		$f^{k}$ 	& ($2^{-k} f^{\log},f^{\arg})$ \\
%JPZ		$h^{k+1}$	& ($2^{-k-1} h^{\log},h^{\arg})$ \\
%JPZ		$c^{k+1}$       & ($2^{-k-1} c^{\log},c^{\arg})$
%JPZ	\end{tabular}
%JPZ	}
\begin{equation} \label{defren}
\left\{ \begin{array}{lll}
f^{k}     & \bydef    & (2^{-k} f^{\log},f^{\arg}) \\
h^{k+1}   & \bydef    & (2^{-k-1} h^{\log},h^{\arg}) \\
\end{array}
\right.
\end{equation}

We also introduce operators:
\begin{eqnarray}
 (x,\alpha) \rentimes{k} (y, \beta) &\bydef& 
 (x + y, \alpha + \beta) \\
 (x,\alpha) \renpow{k}{\lambda} &\bydef& 
 (\lambda x , \lambda \alpha) \\
 \renscal{k}{z} (x,\alpha) &\bydef& 
 (x + 2^{-k} \log |z|, \alpha + \arg z) 
\end{eqnarray}
and
\begin{small}
\begin{eqnarray}
 (x,\alpha) \renplus{k} (y, \beta) &\bydef& 
 \left(
        2^{-k} \log 
 	\left|
	e^{i \alpha + 2^k x}
	+
	e^{i \beta + 2^k y}
	\right|
  ,
  \arg
  \left(
	e^{i \alpha + 2^k x}
	+
	e^{i \beta + 2^k y}
  \right)
  \right) \label{deflast}  
\end{eqnarray} \end{small}\newtext{}{A6\S 1(24) Added parenthesis}
\par
We remark that the purpose of the sub-index $k$ in
the above formulae is to keep track of the degree of the 
renormalization. For operations which do not change,
we omitted the sub-index. \newtext{The operator
$\renscal{k}{z}$ stands for the multiplication of
a renormalized value by a (non-renormalized) constant $z$}{A6\S 1(23)} 
\par
 Also, binary renormalized operations are defined for operands
 with the same renormalization index. Therefore, one should
 first convert $f_i^k$ to $f_i^{k+1}$ before attempting to
 \oldtext{to}\newtext{}{B\S 4(9)}
 `multiply' it with a factor of renormalization index
 of order $k+1$. This conversion will be implicit in the
 formulae below.
\par 
	Equation~(\ref{g2}) becomes:
\begin{small}
\begin{equation}
\label{g666}
{h_i}^{k+1} = 
\left( \renscal{k+1}{ s_{i0}}
\left({f_i}^k \right)
\renpow{k+1}{2} \right)
\renplus{k+1}
\left( \renscal{k+1}{2}
%  {\substack {{\bigrenplus{k+1}}\\
%             \text{\scriptsize $1 \le j \le \min(i, d-i)$} \\}}
\overset{\min(i, d-i)}{\underset{j=1}{\ \ \ \ \bigrenplus{k+1}} } 
\left(
\renscal{k+1}{s_{ij}}
{f_{i-j}}^k
\rentimes{k+1}
{f_{i+j}}^k
\right) \right) \mbox{ ,}
\end{equation} \end{small}
where $s_{ij}= (-1)^{d+i-j}$.
	\newtext{Recall that above, $k$ is a superscript, not an
	exponent.}{A6\S 2(3)}
	The `renormalized operations' above are \oldtext{very}
	\newtext{}{A6\S 2(1)} easy to
	implement in terms of the classical ones. The most 
        delicate being the renormalized sum, so we give here our
	preferred algorithm

\begin{example} \label{ex2} How to compute the `Renormalized sum':
\[
	(c,\gamma) := (a,\alpha) \renplus{k} (b,\beta)
\]\newtext{}{A7(13)}
\begin{tt} 
\begin{tabbing}
\ \ \ \ \= \\
\> If  $a > b$, \= do:  \\
\>		\>$ s = \exp{i\alpha} + \exp{ (i \beta + 2^{k} (b-a)) } $ \\
\>		\>$ c = a + 2^{-k} \ln |s| $ \\
\>		\>$ \gamma = \arg(s) $ \\
\> else  \\
\>		\> $ s  = \exp{i\beta} + \exp{ (i \alpha + 2^{k} (a-b)) }$  \\
\>		\>$ c = b + 2^{-k} \ln |s| $  \\
\>		\>$ \gamma = \arg(s) $  \\
\> endif\\
\end{tabbing}
\end{tt}
\end{example}
\oldtext{$ c = a + 2^{-k} \ln |s| $}
We remark that in the above formula, the complex arithmetic 
operations can be performed in terms of real elementary 
ones. Moreover, in numerical implementations if $k$ is large 
enough, as compared to $\epsilon$, it may be faster to approximate
$(c,\gamma)$ with $(a,\alpha)$ \newtext{or $(b,\beta)$, whichever is
larger}{A7(14)}.

	\medskip
	\par	In order to finish
		the proof of Theorem~\ref{th1}, we
		still need a few remarks about the renormalized Graeffe 
        	\oldtext{operator} 
	                \newtext{iteration, which}{A7(3)}
		we just constructed. 
		Clearly, in order to compute formula~(\ref{g666}),
		we only need memory space of the order $\bigo(d)$
		\newtext{and time of the order $\bigo(d^2)$}{A6\S 1(25)}.
	\par	When $k \rightarrow \infty$, the quantities 
                $|f^k_{i}|$ are all convergent.

	\medskip
	\par
  	If $G_k$ is the renormalized Graeffe iteration, we 
  	define $G_{\infty}$ as the limit of $G_k$ when $k$ goes to infinity.
  	This means that renormalized sums $\renplus{k}$ are replaced by
  	their limit $\renplus{\infty}$, where:
\[
  (a, \alpha) \renplus{\infty} (b, \beta) =
  \left\{
    \begin{array}{ll}
    (a, \alpha) & \text{if } a>b \\
    (b, \beta) & \text{if } a<b \\
    \text{undefined} & \text{if } a=b
    \end{array}
  \right.
\]
\par
	It is easy to see that $G_k \rightarrow G_{\infty}$ almost
	everywhere, \oldtext{and that convergence is locally $C^1$ almost
	everywhere. }\newtext{pointwise and in the $C^1$ topology.
	In order to avoid a rather tedious calculation, we can establish
	this fact from the pointwise $C^1$ convergence almost
	everywhere of
	$(a, \alpha) \renplus{k} (b, \beta)$ to 
	$(a, \alpha) \renplus{\infty} (b, \beta)$ and similarly for
	the other renormalized operations. 
	}{A4\S 2(1-2)}

\newtext{
	We define:
\begin{eqnarray*}
\psi(f_0, \cdots, f_d) 
	&=&
	\left( \log|f_0|, \cdots, \log|f_d| ;
	\arg|f_0|, \cdots, \arg|f_d| \right) \\
\eta(r_0, \cdots, r_d ; a_0, \cdots, a_d) 
	&=&
	\left( \exp {(r_0-r_1)}, \cdots, \exp {(r_{d-1} - r_d)} 
	\right) \\
R(r_0, \cdots, r_d ; a_0, \cdots, a_d) 
	&=&
	\left( \frac{r_0}{2}, \cdots, \frac{r_d}{2}
        ; a_0, \cdots, a_d \right) \\
\end{eqnarray*}
\medskip
\par
	We define $\Phi$ as the function that associates
	to a polynomial $f$ the values $|\zeta_1|, \dots, |\zeta_d|$
        where $\zeta_i$ are the roots of $f$, ordered by decreasing
	modulus. 
\par
	Then, we set
\[
\varphi = \eta^{-1} \circ \Phi \circ \psi^{-1}
\]
\medskip
\par
	With the definitions above, \oldtext{$G^{R,k}$} 
	\newtext{$G_k$}{A3\S 3(6)} is indeed a renormalization
	of the classical Graeffe algorithm to compute $\Phi$: 
}{A6\S 1(26)} 

\begin{proposition} \label{prop1}
	Renormalized Graeffe Iteration is a renormalized iterative
	\newtext{ algorithm (in the sense of Definition~\ref{defren2})  
to compute $\Phi$}{A6\S 1(26)}.
\end{proposition}

	The proof of this proposition will conclude the proof
	of Theorem~\ref{th1}. In order to establish Proposition~\ref{prop1},
	we should check that our algorithm, as described in equation
	(\ref{g666}), satisfies Axioms 1 to 4.
\medskip
\par
	Axiom 1 is verified \newtext{by construction.}{}\oldtext{because of diagram~(\ref{cd}).}
	Axiom 2 \oldtext{is obvious}\newtext{follows}{A6\S 2(2)}
	from the recurrence formula~(\ref{g666}).   
\medskip

\par	The proof that our algorithm satisfies Axiom 3 will require
	a technical lemma. Before stating it, we recall some notation:

\par
   $\text{orb}(f) = \{ f ; G_1(f) ; G_2 \circ G_1 (f) ; \dots 
  ;  G_k \circ \cdots \circ G_1 (f) ; \dots \}$
  \newtext{}{B\S 4(10)} is the orbit of $f$
  under the sequence $(G_k)$. Its closure is denoted by  
  $\overline{ \text{orb}(f) }$.

\oldtext{
   Also, we say that a set is {\em generic} when its complement
   has measure zero.
}\newtext{}{Already defined}
   We will show the following lemma, which implies satisfaction of
	Axiom 3.

{\lemma \label{technical} 
   \oldtext{ Let $f$ be generic 
   and let $\delta > 0$.} \newtext{ For almost every $f$ and 
any $\delta>0$,}{A3\S 1} there exist a 
   compact neighborhood $W \subseteq \mathbb R^l \times \mathbb T^m$
   of $\text{orb} (f)$ and an integer $k_0$
   such that, $G_k$ is a local diffeomorphism $W \rightarrow G_k(W)$,
   and such that for $k_0 \le k \le \infty$, the derivative of $G_k$ is
   bounded by $2+\delta$
}

\begin{proof}[Proof of Lemma~\ref{technical}]
  The proof is divided in several steps.
\begin{itemize}
\item [Step 1:]For any $j \in \mathbb N$, there is a \oldtext{generic}
      \newtext{full-measure set}{A3\S 1} $U_j$ such that
      $G_j \circ \cdots \circ G_1$ is well-defined, and a local
      diffeomorphism in $U_j$. Hence, there is a \oldtext{generic}
      \newtext{full-measure set}{A3\S 1} $U_{j,k}$ such that
\oldtext{      
      $G_k$ 
      is well-defined and a local diffeomorphism near
      $G_j \circ \cdots \circ G_1 (f)$.}
\newtext{$G_k \circ \left( G_j \circ \cdots \circ G_1 \right)$ 
      is defined on $U_{j,k}$ and for all $f \in U_{j,k}$ there is a
      neighborhood $V_f$ of $G_j \circ \cdots \circ G_1 (f)$ such that
      $G_k$ is a diffeomorphism $V_f \rightarrow G_k(V_f)$.}
      {A3\S 2(2)}
\item [Step 2:]Let $U_{\infty}$ be the set of all $f$ such that 
      \oldtext{
      $G_k(f)$ is
      a local diffeomorphism in $(\pi ^{-1} \circ \varphi) (f)$, for
      $k$ large enough.}
      \newtext{$G_k$ is a diffeomorphism near 
      $(\pi ^{-1} \circ \varphi) (f)$ for all values of $k$ that
      are large enough.}{A4\S 1(9)}
      Then $U_\infty$ contains the set of complex 
      polynomials without roots of the same modulus. Hence $U_\infty$ 
      \oldtext{is generic.} \newtext{has full measure.}{A3\S 1}
\item [Step 3:]Let $U = U_\infty \cap \left( \bigcap_j U_j \right)  
      \cap \left( \bigcap_{jk} U_{jk} \right)$. Then $U$ 
      \oldtext{is generic.} \newtext{has full measure}{A3\S 1}.
      Moreover, let $f \in U$. Then $G_k$ is a local
      diffeomorphism with derivative \newtext{of norm}
      {This is more precise.} $< 2+\delta / 2$ in an open
      neighborhood $V_k(g)$ of every $g \in \overline{\text{orb} f}$,
      for all $k \ge k_0 (g)$. 
      \newtext{Indeed, if we write $G_k = G_{\infty} + 
      \left( G_k - G_{\infty} \right)$, we can make the
      $C^1$ norm of the second term arbitrarily small, namely
      less than $\delta /2$. We know that the norm of the derivative
      of $G_{\infty}$ is precisely 2, hence the bound $2+\delta / 2$
      }{A4\S 2(1)}
      In the particular case $\delta = 
      \infty$ we can set $k_0 = 1$.
\item [Step 4:]Since $G_k \rightarrow G_{\infty}$ \oldtext{in a locally
	$C^1$ sense}\newtext{pointwise in the $C^1$ topology and }{A4\S 2(2)}
	for $g$ almost everywhere, we can assume 
	that
      $\bigcap_{k \ge k_0} V_k(g)$ contains an open ball $V(g)$ of
      center $g$, where $G_k$ is a local diffeomorphism with 
      derivative bounded by $2+\delta / 2$. 
\item [Step 5:]Since $\overline {\text{orb}(f)}$ is compact, the union
      $\bigcup_{g \in \overline{\text{orb}(f)}} V(g)$ has a finite 
      sub-cover $\bigcup_{g \in \Gamma} V(g)$, and we set $W =
      \bigcup_{g \in \Gamma} \overline{V(g)}$. Then we set $k_0
      = \max_{g \in \Gamma} k_0(g)$, and we obtain that for any
      $k \ge k_0$, $G_k$ is a local diffeomorphism in $W$, with
      derivative of norm bounded by $2 + \delta$. 
\end{itemize}
\end{proof}

\medskip
\par
	The proof that our algorithm satisfies Axiom 4 is divided in two parts
	\newtext{, dealing (respectively) with small and large values of $k$}
	{For clarity.}.
\oldtext{
\begin{lemma*}
%	\label{lp1}
	For a.e. $g$ given, $|G_{k,\epsilon}(g) - G_{k}(g) | < B \epsilon$,
	where $B$ depends on $g$ and $k$, and $\epsilon$ is small
	enough.
\end{lemma*}

\begin{lemma*} 
% \marginpar{This Lemma and its proof need to be rewritten.}
%	\label{lp2}
	For a.e. $g$ given, there is $k_0$ such that for any 
	$k > k_0$
\[
\dist \left(G_{k,\epsilon}(g) , G_{k}(g) \right) < B \epsilon
\]
	where $B$ depends on $\varphi(g)$ only.
\end{lemma*}
}
\newtext{
\begin{lemma}
	\label{lp1}
	For almost every $f$ and for all $k$,
	there exist an open neighborhood
	$U$ containing $f$ and $B>0$ such that for all
	$g \in U$ and for all $\epsilon$ small enough
\[
\dist \left(G_{k,\epsilon}(g) - G_{k}(g) \right) < B \epsilon \mbox{ .} 
\]
\end{lemma}

\begin{lemma} 
	\label{lp2}
	For almost every $f$, there exist $k_0$, $B>0$, 
	and an open neighborhood
	$U$ containing $f$, such that for all
	$g \in U$ and $k \ge k_0$, for all $\epsilon$ small enough,
\oldtext{
\[
|G_{k,\epsilon}(g) - G_{k}(g) | < B \epsilon
\]
	For a.e. $g$ given, there is $k_0$ such that for any 
	$k > k_0$
}
\[
\dist \left(G_{k,\epsilon}(g) , G_{k}(g) \right) < B \epsilon
\]
\end{lemma}
}{Restated for clarity and also A3\S 3(7).}

\medskip
\par	\oldtext{Since $G_k$ is an almost everywhere local diffeomorphism, 
        and hence the preimage
	of a zero measure set has zero measure,} Lemmas~\ref{lp1} and
	\ref{lp2} together imply that for almost all \oldtext{$g$}
	\newtext{$f$, there is a neighborhood $U$ containing $f$ 
	such that for all $g \in U$}{},
\[
\dist \left( G_{k,\epsilon}(g) , G_{k}(g) \right) < B \epsilon
\]
	independently of $k$.
\par

	\oldtext{
	This result can be extended to an open neighborhood
	of any $g$ generic.}
	Since \oldtext{(generically)} $\text{orb}(f)$ is defined
        \newtext{for almost every $f$}{A3\S 1} 
	and admits a compact neighborhood $V$, we can select a finite
        subcover of those neighborhoods, and hence find a finite $B$
        valid for all $g \in V$. Therefore, our algorithm satisfies
        Axiom 4. 
\medskip
\par	
\begin{proof}[Proof of Lemma~\ref{lp1}:]
	Our iteration $G_k$ can be written in terms of the following
	real operations: $+$, $-$, $\cos$, $\sin$, $\arctan$, 
	multiplication by $2^k$, by $2^{-k}$, absolute value,
	$\exp$, $\log$.
\par	The set of inputs such that a `log of zero' or an `absolute
	value of zero' occurs has zero measure. \oldtext{Therefore, if $g$ is
	generic,}
	\newtext{Therefore, for almost every $f$,}{A3\S 1}	
	$G_k(f)$ is computed by a composition of analytic
	functions. Also, \oldtext{if $g$ is generic}
	\newtext{for almost every $f$,}{A3\S 1} 
	none of the output values
	is zero.
\par	Therefore, for every intermediate quantity $x_l$, the derivative
	of any output $y_m$ with respect to $x_l$ is finite (say $\le D_{lm}$).
\oldtext{
	Therefore, a relative perturbation of $\epsilon$ in $x$ leads
	to a perturbation of $xD\epsilon$ in $y$. Thus, we set
	$B = \sum D_{lm}x_l$ and we are done. 
	}
\newtext{
	Therefore, a relative perturbation of $\epsilon$ in $x_l$ leads
	to a perturbation of size $|x_l| D_{lm}\epsilon$ in $y_m$. Thus, we set
	$B = \frac{1}{2} \sum_{l,m} D_{lm}|x_l|$. By continuity,
	a perturbation of $\epsilon$ in each intermediate value $x_l$
	leads to a perturbation smaller than $B \epsilon$, for input
	$g$ in a certain neighborhood of $f$. 
	}{A6\S 1(27)}
\end{proof}
\medskip
\par	
\begin{proof}[Proof of Lemma~\ref{lp2}:]
\oldtext{
	We claim that there is a full-measure set $W$ such that 
	if $\varphi(g) \in W$, then there are $k_0$ and $B$ such 
	that 
\[
\dist \left( G_{k,\epsilon}(g) , G_{k}(g) \right) < B \epsilon
\]
	for $k > k_0$. Since $\varphi ^{-1}$ maps zero
	measure sets into zero measure sets, Lemma~\ref{lp2} will
	be true in the full-measure set 
	$\varphi^{-1}(W) \cap W'$. where $W'$ is the set of
	$g$ such that $G_k \circ \dots \circ G_1 \rightarrow \varphi(g)$.
	According to Axiom 1, the set $W'$ has full measure.
}
\newtext{
	Let $W$ be the set of all $f$ such that $G_{\infty}(f)$ is
	well-defined. Recall that
\[
	\begin{array}{lccll}
	\renplus{\infty} : & (a,\alpha), (b,\beta) &
	\mapsto & (a,\alpha) & \text{if  } a>b \\
	&&& (b,\beta) & \text{if  } b>a 
	\end{array}	
\]
	and that the operator $\renplus{\infty}$ is not defined 
	for $a=b$. Therefore, $W$ is open and has full measure.}{Clarity.}
\oldtext{
\par	In order to prove our claim, we consider the `limit' iteration
	$G_{\infty}$, where $\renplus{k}$ was replaced by 
	$\renplus{\infty}$:
\[
	\begin{array}{lccll}
	\renplus{\infty} : & (a,\alpha), (b,\beta) &
	\mapsto & (a,\alpha) & \text{if  } a>b \\
	&&& (b,\beta) & \text{if  } b>a 
	\end{array}	
\]
\par	The operator $\renplus{\infty}$ is not defined for $a=b$.
	In the same way, we replace all the other renormalized operators by
	their limit values.
	Wherever defined, $G_{\infty}$ fixes scalar quantities, and
	doubles angles.
\par	Let $W$ be the domain of definition of $G_{\infty}$. Clearly,
	$W$ is open and has full measure.}
\newtext{
	Let $f \in W$ and $U$ be a small connected
	neighborhood containing $f$.
	Let $g \in U$, then
	by taking $U$ small enough and $k$ large enough, we
	can guarantee that $G_k(g)$ and $G_{k,\epsilon}(g)$
	are well-defined. Since $U$ is connected, all the
	branching outcomes in the computation of $G_k(g)$ and
	$G_{k,\epsilon}$ are
	the same, hence $G_k$ restricted to $U$ is a composition
	of locally analytic functions. We can assume without loss
	of generality that all derivatives are bounded, hence
	there is a constant $B$ such that 
\[
	\dist \left(G_{k,\epsilon}(g) , G_{k}(g) \right) < B \epsilon
\]
	for $\epsilon$ small enough, but still independent
	of the choice of $g$.
	}{Fixes A3\S 3(7). The following remarks are no more relevant:
        A3\S 3(8), A6\S 1(28), A6\S 1(30), A3\S 2(3)}
	
\oldtext{	
	For $y \in W$, there is $k$
	such that every operation $\renplus{k}$ can be replaced by
	$\renplus{\infty}$, with relative error $\epsilon$. Indeed,
	for $k$ large enough,
\[
	s = \left| 1 + e^{i(\beta - \alpha) + 2^k (b-a)} \right| \ , \ a > b
\]
	is bounded away from zero, so that 
	$2^{-k} \log s$ can be
	made smaller than $|a| \epsilon$. The error in the angle may be
	estimated the same way. 
\par	This last property holds in an open neighborhood $U$ of $y$,
	$U \subset W$. Once we take $k > k_0$, with $k_0$ large enough,
        $(G_k \circ \dots \circ G_1) (f) \in U$. Since all of our
	operators $\renplus{k}$ and $\rentimes{k}$ are 
	locally Lipschitz in
	$W$, there is $B$ s.t. 
        for $\epsilon$ small enough, for all $z \in U$,
$\| G_{k \epsilon}(z) - G_k(z) \|$
is uniformly bounded by $B \epsilon$}
\end{proof}	
\medskip
\par
\oldtext{
\begin{proof}[End of the proof of Theorem~\ref{th1}:]
	We define:
\begin{eqnarray*}
\psi(f_0, \cdots, f_d) 
	&=&
	\left( \log|f_0|, \cdots, \log|f_d| ;
	\arg|f_0|, \cdots, \arg|f_d| \right) \\
\eta(r_0, \cdots, r_d ; a_0, \cdots, a_d) 
	&=&
	\left( \exp {(r_0-r_1)}, \cdots, \exp {(r_{d-1} - r_d)} 
	\right) \\
R(r_0, \cdots, r_d ; a_0, \cdots, a_d) 
	&=&
	\left( \frac{r_0}{2}, \cdots, \frac{r_d}{2}
        ; a_0, \cdots, a_d \right) \\
\end{eqnarray*}
\medskip
\par
	We define $\Phi$ as the function that associates
	to a polynomial $f$ the values $|\zeta_1|, \dots, |zeta_d|$
        where $\zeta_i$ are the roots of $f$, ordered by decreasing
	modulus. 
\par
	Then we set
\[
\varphi = \eta^{-1} \circ \Phi \circ \psi^{-1}
\]
 
\medskip
\par
	With the definitions above, $G^{R,k}$ 
	is indeed a renormalization
	of the classical Graeffe algorithm to compute $\Phi$. 
\end{proof}
}
\newtext{This concludes the proof of Proposition~\ref{prop1} and
hence of Theorem~\ref{th1}}{A6\S 1(26)}

\section{Probability of Success}
\label{Sec5}

\oldtext{MOST OF THIS SECTION WAS REWRITTEN FROM SCRATCH}

Through this section, $\|.\|_d$ will denote Weyl's unitary
invariant norm (\cite{WEYL} ~III--7, pp~137--140) 
in the space $\mathcal P_d$ of complex 
polynomials of degree at most $d$ (See~\cite{BCSS}). 
If $f(x) = \sum_{i=0}^d f_i x^i$,
then
\[
\| f \|_d = \sqrt{ \sum_{i=0}^d \frac{|f_i|}{\binomial{d}{i}}}
\]
\par
This norm is invariant under the following action of the group
$U(2)$ of unitary $2\times 2$ matrices: if 
$\varphi = \left[ \begin{array}{cc} \alpha & \beta \\ \gamma & \delta
\end{array} \right]$ is unitary, define
\[
f^{\varphi} (y) = (\gamma y + \delta)^d f\left( \frac{\alpha y + \beta}
{\gamma y + \delta} \right)
\]
Thus,  $\| f^{\varphi} \|_d = \| f \|$. (See Ch.~12, Th.~1 of \cite{BCSS}.)

This is in some sense the most `natural' norm in the space of
all polynomials.  \newtext{More information about that norm, its associated
probability distribution and its applications can be found in}{A5(14)}
\cite{BD,DEDIEU2,DEDIEU3,KOSTLAN,TCS,SPLITTING,BEZI,BEZII,BEZIII, BEZIV,
BEZV}.

Let $\mathbb P \mathcal P_d$ be the projectivization of normed complex vector
space $(\mathcal P_d, \|.\|_d)$. We can define the `sine' distance in
$\mathbb P \mathcal P_d$ by:
\[
d_{\mathbb P} (f,g) = \min_{\lambda \in \mathbb C} \frac{\|f - \lambda g\|}
{\|f\|} 
=
\sin \varrho(f,g)
\]
where $\varrho$ is the usual Riemann metric in $\mathbb P \mathcal P_d$.
Also, there is a natural volume form in $\mathbb P \mathcal P_d$.
Normal invariant distributed polynomials in $\mathcal P_d$ correspond
to uniformly distributed `polynomials' in $\mathbb P \mathcal P_d$.

Let $f$ be a random polynomial (in any of the two equivalent senses above).
Then we may consider its roots $\zeta_1, \cdots, \zeta_d$ as random variables.
The joint distribution of $\zeta_1, \cdots, \zeta_d$ was studied by Kostlan in
~\cite{KOSTLAN}. However, the result below will rely on elementary
estimates:

\begin{lemma} \label{lem6}
  Let $f$ be random in the sense above.
  The probability that 
\[
\min_{|\zeta_i|>|\zeta_j|} \frac{|\zeta_i|}{|\zeta_j|} > 1+\epsilon
\]
  is larger than $1-M \epsilon$, where $M$ is a positive constant depending
  on $d$.
\end{lemma}

  This result is also true if we chose $f$ random with respect to any 
  other probability distribution \oldtext{
  
  that is absolutely continuous
  with respect to the volume form of $\mathbb P \mathcal P_d$.}
  \newtext{with bounded Radon-Nikodym derivative with respect to
  the volume form in $\mathbb P \mathcal P_d$.}{Precise definition.}
  The constant $M$ will have to be multiplied by the maximum of the
  Radon-Nikodym derivative of the new distribution with respect to
  the volume form.

Lemma~\ref{lem6} will be a consequence of the `Condition Number 
Theorem' below. Let
\[
\rho(f) = \min_{|\zeta_i| < |\zeta_j|} 1 - \frac{|\zeta_i|}{|\zeta_j|}
\]

We will interpret $\rho(f)^{-1}$ as a condition number. Let
$\Sigma_G$ be the locus of ill-posed problems, i.e., the set of
polynomials such that $|\zeta_i| = |\zeta_j|$ for some $i \ne j$. Then,
\begin{theorem}[Condition Number Theorem for Graeffe Iteration]\label{cnt}
\[
\rho(f) \ge \frac{ d_{\mathbb P} (f, \Sigma_G) }{\sqrt{d}}
\]
\end{theorem}

Therefore, the probability that $\rho(f)>\epsilon$ is no
less than $1 - \vol V_{\epsilon \sqrt{d}} \Sigma_G$ where 
$V_{\epsilon \sqrt{d}}$ denotes an $\epsilon \sqrt{d}$-neighborhood. 
This volume is of $O(\epsilon \sqrt{d} ) \le M \epsilon$.

\begin{lemma}\label {norms1}
	Let $f(x) = (x-\zeta_1) g(x) \in \mathcal P_d$, where
	$g \in \mathcal P_{d-1}$. Then
\[
\| g \|_{d-1} \le \sqrt{ \frac{d}{1+|\zeta_1|^2}} \| f \|_d
\]
\end{lemma}

\begin{proof}[Proof of Lemma~\ref{norms1}]

We start with the easy case and assume that $\zeta_1=0$. Set
$g(x) = \sum_{i=0}^{d-1} g_i x^i$. Then, $f(x) = \sum_{i=1}^{d} g_{i-1} x^i$.
Now,
\[
\|g\|_{d-1}^2 
=
\sum _{i=0}^{d-1} \frac{|g_i|^2}{\binomial{d-1}{i}}
=
\sum _{i=1}^{d} \frac{|g_{i-1}|^2}{\binomial{d-1}{i-1}}
=
\sum _{i=1}^{d} \frac{d}{i} \frac{|g_{i-1}|^2}{\binomial{d}{i}}
\le
d \|f\|_d ^2
\]
and hence $\|g\|_{d-1} \le \sqrt{d} \|f\|_d$.

For the general case, we will use $U(2)$-invariance of
$\|.\|_d$ and $\|.\|_{d-1}$. Let $\varphi$ be a convenient
unitary matrix:
\[
\varphi = \frac{1}{\sqrt{1+|\zeta_1|^2}}
\left[
\begin{array}{cc} 1 & \zeta_1 \\ - \bar{\zeta_1} & 1 \end{array} \right]
\]

Set $f^{\varphi} = f \circ \varphi$ and $g^{\varphi} = g \circ \varphi$.
The choice of $\varphi$ has the particularity that $f^{\varphi}(0) =
f(\zeta_1) = 0$. We can compute $f^{\varphi}$ in terms of $g^{\varphi}$:
\[
f^{\varphi}(y) = \sqrt{1 + |\zeta_1|^2 } \ y \ g^{\varphi} (y)
\]

Using the easy case,
\[ 
\|  \sqrt{1 + |\zeta_1|^2 } g^{\varphi} \|_{d-1}
\le
\sqrt{d} \|f^{\varphi}\|_d
\]

By $U(2)$-invariance,

\begin{eqnarray*}
\| g \|_{d-1} &=& \| g^{\varphi} \| _{d-1} \\
&=& \frac{1}{ \sqrt{1 + |\zeta_1|^2 }}  \| \sqrt{1 + |\zeta_1|^2 } g^{\varphi} \| _{d-1} \\
&\le& \frac{\sqrt{d}}{ \sqrt{1 + |\zeta_1|^2 }}  \| f^{\varphi} \|_d \\
&=&
\frac{\sqrt{d}}{ \sqrt{1 + |\zeta_1|^2 }}  \| f \|_d
\end{eqnarray*}
\end{proof}

\begin{lemma} \label{norms2}
 Let $g \in \mathcal P_{d-1}$. Then $\|g\|_d \le  \|g\|_{d-1}$.
\end{lemma}

\begin{proof}[Proof of Lemma~\ref{norms2}]
\[
\|g\|_d^2 = 
\sum_{i=0}^{d-1} \frac{|g_i|^2}{\binomial{d}{i}}
=
\sum_{i=0}^{d-1} \frac{d-i}{d} \frac{|g_i|^2}{\binomial{d-1}{i}}
\le
\|g\|_{d-1}^2 
\]
\end{proof}

Putting Lemma~\ref{norms1} and Lemma~\ref{norms2} together,

\begin{lemma}\label {norms3}
	Let $f(x) = (x-\zeta_1) g(x) \in \mathcal P_d$, where
	$g \in \mathcal P_{d-1}$. Then
\[
\| g \|_{d} \le \frac{\sqrt{d}}{\sqrt{1+|\zeta_1|^2}} \| f \|_d
\]
\end{lemma}

\begin{proof}[Proof of Theorem~\ref{cnt}]

Let $f(x) = (x-\zeta_1) (x-\zeta_2) \cdots (x-\zeta_d)$ and order the $\zeta_i$'s
such that
\[
\rho(f) = \min_{|\zeta_i| < |\zeta_j|} 1- \frac{|\zeta_i|}{|\zeta_j|} 
= 1 - \frac{|\zeta_2|}{|\zeta_1|} 
\]

Define $h(z) = (x - \zeta_1 (1-\rho)) (x-\zeta_2) \cdots (x-\zeta_d) \in \Sigma_G$.
Then
\begin{eqnarray*}
\|f - h\|_d
&=&
\| \zeta_1 \rho (x-\zeta_2) \cdots (x-\zeta_d) \|_d \\
&=& |\zeta_1| \rho \| (x-\zeta_2) \cdots (x-\zeta_d) \|_d\\
&\le&
|\zeta_1| \rho \frac{\sqrt{d}}{\sqrt{1+|\zeta_1|^2}} \|f\|_d
\end{eqnarray*}

Hence,

\[
\frac{\|f-h\|_d}{\|f\|_d} \le \frac{|\zeta_1|}{\sqrt{1+|\zeta_1|^2}} 
\rho \sqrt{d} \le \rho \sqrt{d} 
\]

Hence,
\[
d_{\mathbb P} (f, \Sigma_G) \le \rho \sqrt{d}
\]
\end{proof}

\begin{proof}[Proof of Theorem~\ref{th2}]

We set $\delta = M \epsilon$, where $M$ is the constant such that
the volume of an $\epsilon \sqrt{d}$ neighborhood of $\Sigma_G$ is less
than $M \epsilon$. With probability
larger than $1-\delta$,
\[ 
\max_{|\zeta_{i}| < |\zeta_{j}|} \frac{|\zeta_{i}|}{|\zeta_{j}|}
> 1 - \epsilon 
\mbox{  .} \]
\par
We now use the fact that, for any $N > (\log 2) / \epsilon$, 
we have $(1-\epsilon)^N < 1/2$. We set $k_1 = 1 + \lceil \log_2 \epsilon^{-1}
\rceil$ and using $N=2^{k_1}$ in the previous formula, we obtain:
\[
\max_{|\zeta_{i}| < |\zeta_{j}|}
\frac{|\zeta_{i}|^{2^{k_1}}}{|\zeta_{j}|^{2^{k_1}}}
<
(1-\epsilon)^{2^{k_1}}
< \frac{1}{2}
\mbox{  .} \]
\par
An extra $1+\log_2 b$ iterations ensures that 
\[
\max_{|\zeta_{i}| < |\zeta_{j}|}
\frac{|\zeta_{i}|^{2^{k}}}{|\zeta_{j}|^{2^{k}}}
<
2^{-1-b}
\mbox{ ,} \]
provided that 
\[ k = k_1 + 1 + \log_2 b \mbox{ .} \]
We claim that we have in this case  that the sum
\[ \sum_{i_1 <  \dots < i_r} \zeta_{i_1}^{2^{k}}\dots \zeta_{i_r}^{2^{k}}= 
\zeta_{1}^{2^{k}} \dots \zeta_{r}^{2^{k}} 
(1+\varepsilon) \mbox{ ,} \] 
where  
\[| \varepsilon |< 2^{-b} \mbox{ .} \]
Indeed, assume we have reordered the roots so that 
$|\zeta_1| >  \dots >  |\zeta_{d}| $.
It then follows that
\[ | \zeta_{i_1}^{2^k} \dots \zeta_{i_r}^{2^k} |  < |\zeta_{1}^{2^k} \dots \zeta_{r}^{2^k}| 
2^{-(1+b)(i_1+\dots + i_{r}- r (r+1)/2)} \mbox{ .} \]

\oldtext{
Now, let 
\[ {\cal  S}_{k} = \{(i_1,\dots,i_n) 
\ \  |  \ \ i_1  + \dots + i_n - n(n+1)/2 =k \} \mbox{ .}  \]
Note that ${\cal  S}_{0}=\{(1,\dots,n)\}$, and
${\cal  S}_{1}=\{(1,\dots,n-1,n+1)\}$. Also,
${\cal  S}_{2}=\{(1,\dots,n-1,n+2) ; (1,\dots,n-2,n,n+1)\}$

But, $\#{\cal S}_{i+1} \le d \# {\cal S}_i$, and hence
\[ \sum_{i\ge 1} 2^{-b i}\#{\cal S}_i \le 2^{-b} (1+2 2^{-b} + 
2 d 2^{-b} + 2 d^2 2^{-2b}\dots ) 
\le 2^{-b+1}\mbox{ ,} \]
under the assumption that 
$d \le 2^b$.
This concludes the proof of Theorem~\ref{th2}.
}
\newtext{
Now, for $r$ fixed, let 
\[ {\cal  S}_{k} = \{(i_1,\dots,i_r) 
\ \  |  \ i_1 < \cdots < i_r \text{ and }
i_1  + \dots i_r - r(r+1)/2 =k \} \mbox{ .}  \]
Note that ${\cal  S}_{0}=\{(1,\dots,r)\}$, and
${\cal  S}_{1}=\{(1,\dots,r-1,r+1)\}$. Also,
${\cal  S}_{2}=\{(1,\dots,r-1,r+2) ; (1,\dots,r-2,r,r+1)\}$
\par
In general, every multi-index $i_1 < \cdots < i_r$ may be
obtained by starting from $1<2< \cdots <r$ and increasing
one of the indices, in such a way not two indices are equal.
Then ${\cal  S}_{k}$ is the set of multi-indices obtained
after $k$ steps. Therefore, we may bound
$\#{\cal S}_{k+1} \le d \ \# {\cal S}_k$ to get
\[ \sum_{k\ge 1} 2^{-(b+1) k}\ \#{\cal S}_k \le 2^{-b-1} (1+2^{-b} + 
2^{-b}d +  2^{-2b-1}d^2 \dots ) 
\le 2^{-b}\mbox{ ,} \]
under the assumption that 
$d \le 2^b$.
This concludes the proof of Theorem~\ref{th2}.
}{A4\S 2(3-4)}
\end{proof}

\section{Newton Diagram Revisited}
\label{Sec4}
 
\begin{figure} 
\vspace{8cm} 
\includegraphics{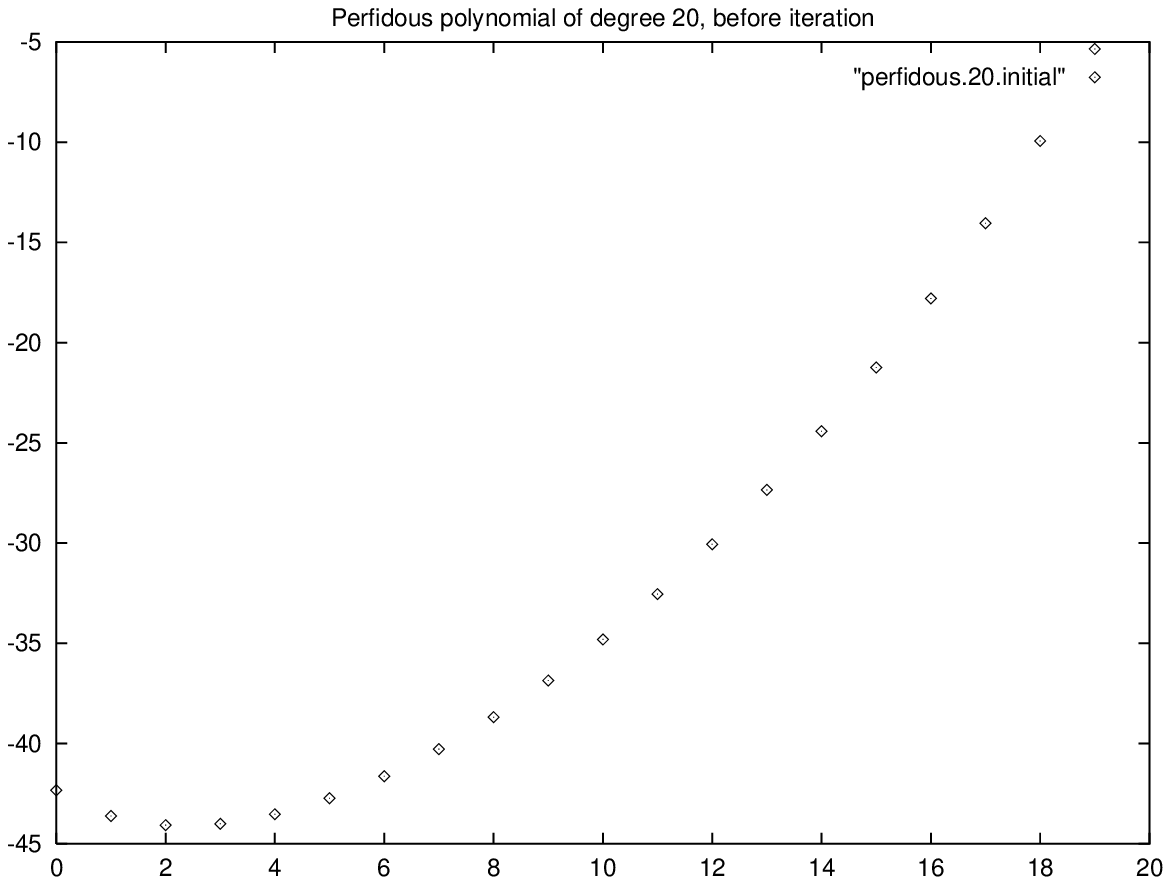}
\includegraphics{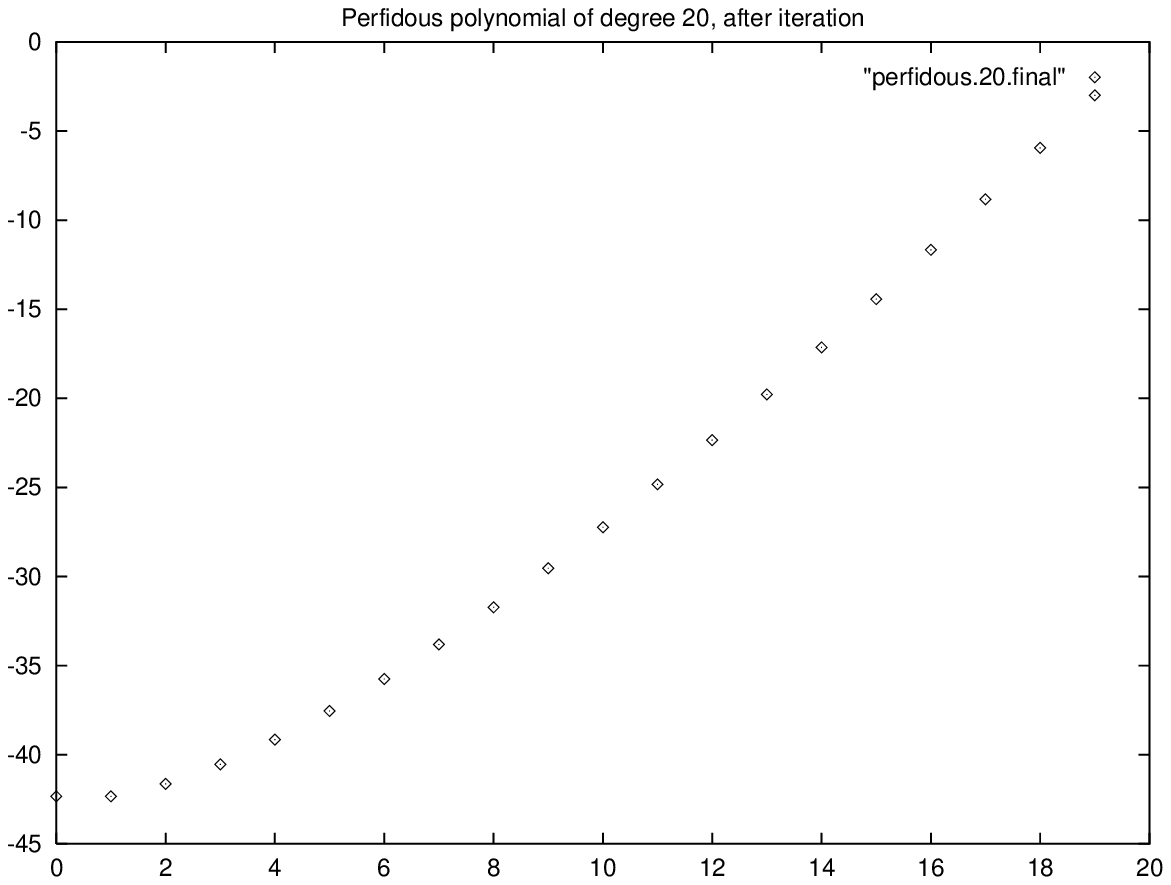} 
\includegraphics{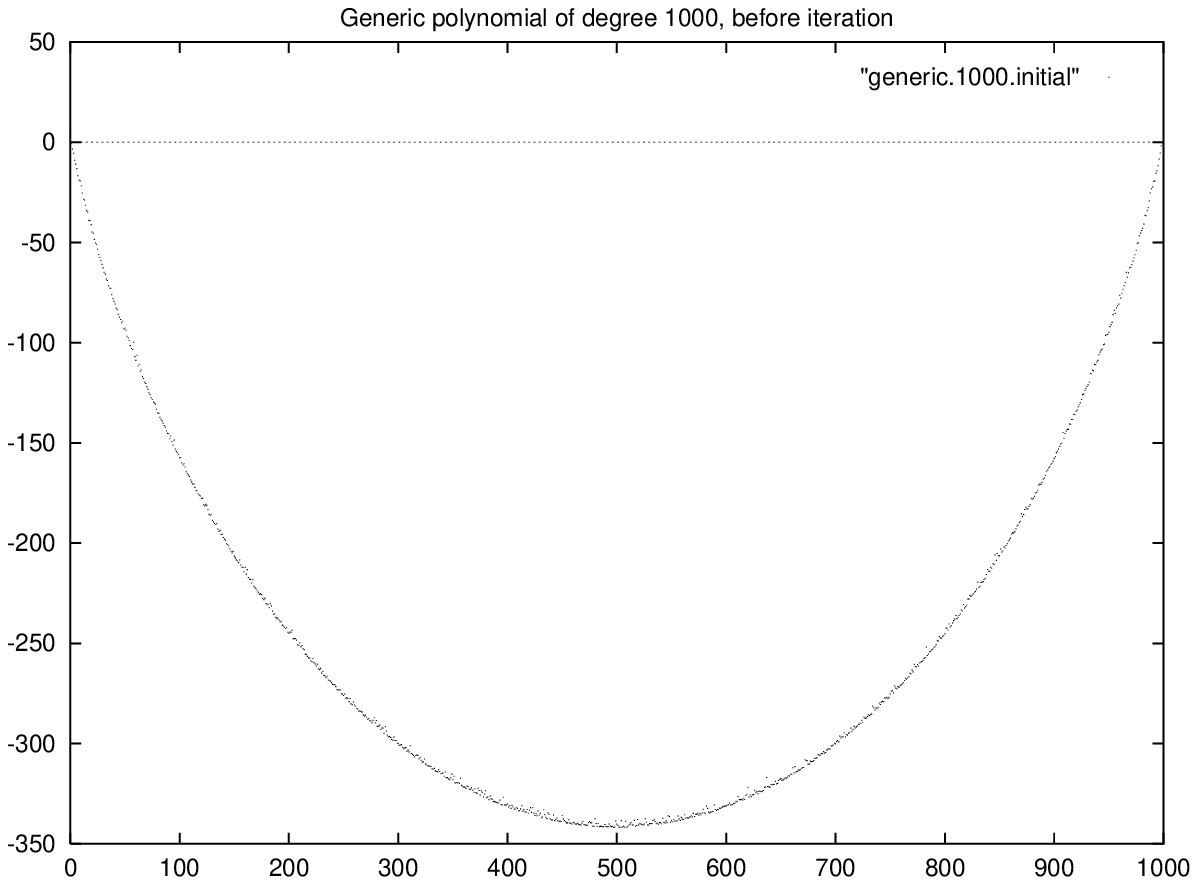}
\includegraphics{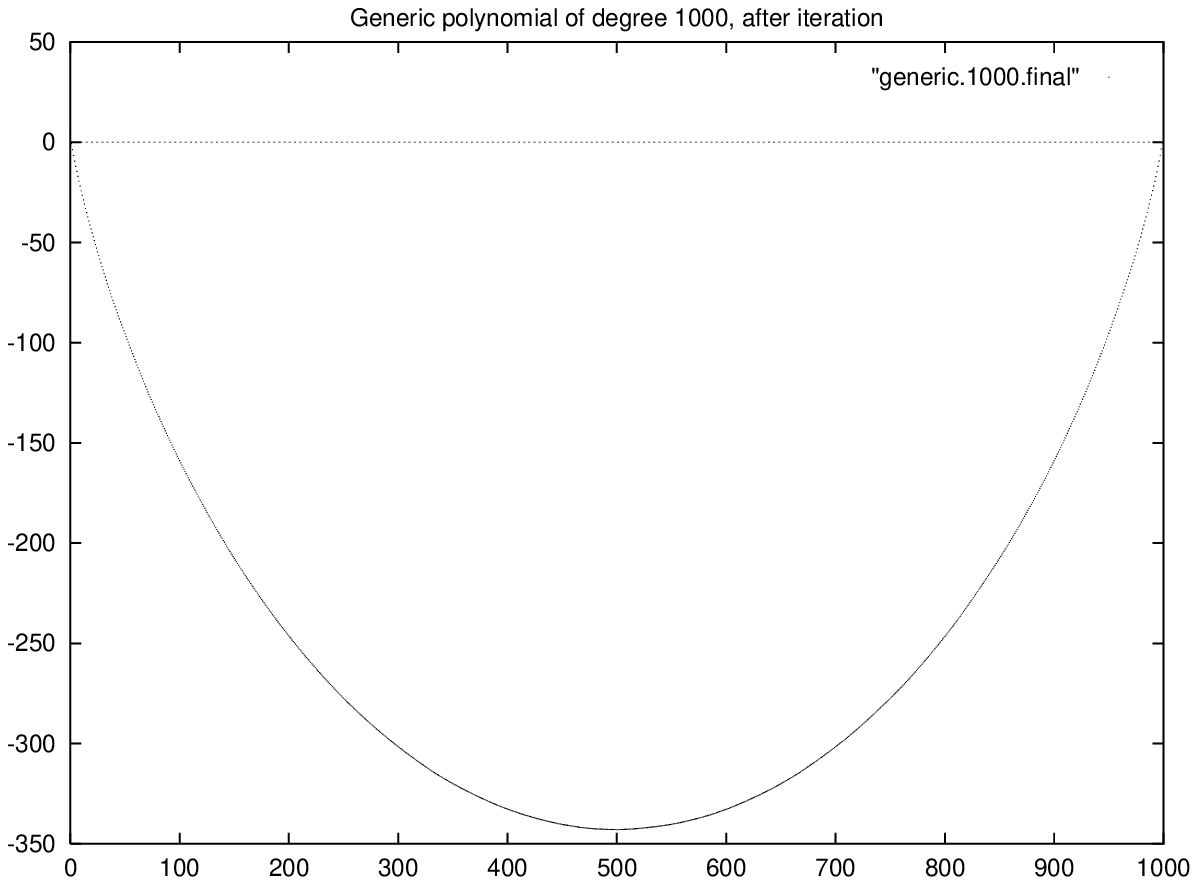}
\ \\
\caption[ ]{\label{newt}Examples of Newton diagrams}
 \end{figure}

	\par	Following Ostrowskii~\cite{OSTROWSKII}, we introduce Newton 
                diagram, but now with a log scale.
		The Newton diagram of a polynomial $f$ is the graph
		of the piecewise linear function defined by $g(i) =
		-\log |f_i|$. 
 (Our construction differs from Ostrowskii's point of view. He used
  to define the Newton diagram as the convex hull of that set of points).
                It makes sense to renormalize Newton
		diagram, as we did with polynomials $G^k f$.
	\par	Therefore, we may speak of a Newton diagram {\em on
		the limit}. Assume (as we did through this point) that
		all roots of $f$ have different moduli. Then the derivative 
		of this diagram on the 
		piecewise linear part over $[i,i+1]$ will correspond to  
		$\log |\zeta_{d-i}|$, where $\zeta_{i}$ is the $i$-th
		largest root of $f$ in modulus. It follows that the 
		Newton diagram
		will be the graph of a convex function.
 	\medskip
	\par	We should consider now the general case. If $f$ has
		several roots of the same moduli, some of the ratios
		$f_i/f_{i+1}$ will converge to the ratios of the coefficients 
		of the factor of $f$ containing those roots. The case of 
		three roots of same moduli is illustrative:
	\par    Let $\zeta_1$, \dots, $\zeta_d$ be arranged by non-increasing
		moduli, and assume $|\zeta_i| = |\zeta_{i+1}| = |\zeta_{i+2}|$.
		Then, 
\begin{eqnarray}
		f_{d-i}   &=& \zeta_1 \dots \zeta_{i-1}
			\left(\zeta_i + \zeta_{i+1} + \zeta_{i+2}\right)
				+ \dots\\
		f_{d-i-1} &=& \zeta_1 \dots \zeta_{i-1} 
			\left( \zeta_i \zeta_{i+1} + 
				+ \zeta_{i} \zeta_{i+2}
				+ \zeta_{i+1} \zeta_{i+2}
				\right)
				+ \dots\\
		f_{d-i-2} &=& \zeta_1 \dots \zeta_{i+2} + \dots 
\end{eqnarray}
	\par	Hence, {\em on the limit}, we approximate:
	\begin{eqnarray}
		\frac{f_{d-i-2}}{f_{d-i-1}}
		&\simeq&
		\frac 
				{\zeta_i \zeta_{i+1} \zeta_{i+2}} 
				{ \zeta_i \zeta_{i+1} + 
				+ \zeta_{i} \zeta_{i+2}
				+ \zeta_{i+1} \zeta_{i+2}
				}
				\\
		\frac{f_{d-i-1}}{f_{d-i}}
		&\simeq&
		\frac  	
			{ \zeta_i \zeta_{i+1} + 
				+ \zeta_{i} \zeta_{i+2}
				+ \zeta_{i+1} \zeta_{i+2}
				} 
			{ \zeta_i + \zeta_{i+1} + \zeta_{i+2}}
			\\
		\frac{f_{d-i}}{f_{d-i+1}}
		&\simeq&
		\zeta_i + \zeta_{i+1} + \zeta_{i+2}
	\end{eqnarray} 
	\par	The moduli of those roots is therefore given by:
	\begin{equation}
		\sqrt[3]{
		\left| \frac{f_{d-i-2}}{f_{d-i-1}} \right|
		\left| \frac{f_{d-i-1}}{f_{d-i}} \right|
		\left| \frac{f_{d-i}}{f_{d-i+1}} \right|
		} 
	\end{equation}		
	\par	The same formula extends for factors of any degree,
		provided they have roots in a circle and all
		other roots are far away from this circle.
	\par	It is useful to have a decision criterion for the
		existence of factors of degree greater than 1. 
		We will do that for degree-2 factors, since this
		is the interesting case for \oldtext{generic}
		\newtext{}{A3\S 1} real polynomials.
		See Ostrowskii~\cite{OSTROWSKII} for more results.
	\par	Clearly, it is enough to consider the case of a 
		polynomial $f$ of degree 2:
\begin{equation}
		f(x) = f_2 x^2 + f_1 x + f_0 
		= x^2 + (- \zeta_1 - \zeta_2) x + \zeta_1 \zeta_2
\end{equation} 
	\par	In case $|\zeta_1| \gg |\zeta_2|$, we have: 
\begin{equation}
\label{r1}
		-\log \left( \frac{|f_2|}{|f_1|} \right) 
		+\log \left( \frac{|f_1|}{|f_0|} \right) \ge 
		\log \frac{|\zeta_1|}{|\zeta_2|} \gg 0 
\end{equation}
	\par	In case $R = |\zeta_1| = |\zeta_2|$, we can bound:
\begin{equation}
\label{r2}
		-\log \left( \frac{|f_2|}{|f_1|} \right) 
		+\log \left( \frac{|f_1|}{|f_0|} \right) < \log 4 
\end{equation}
	\par	However, we should look at the renormalized Newton
		diagram. In that case, we should divide equations
		(\ref{r1}) and (\ref{r2}) by $2^k$
		\newtext{, and expect that if $|\zeta_1|\ne |\zeta_2|$
		then the following holds: }{} 
\begin{equation}
\label{r3}
		-2^{-k} \log \left( \frac{|f_2|}{|f_1|} \right) 
		+2^{-k}\log \left( \frac{|f_1|}{|f_0|} \right) \ge \sigma 
\end{equation}
		where $\sigma$ is an a priori bound on root
		\newtext{moduli}{} separation.
		It may be obtained from a probabilistic analysis
		(Lemma~\ref{lem6}), or from any other a priori  
                knowledge on the polynomial; 
		\newtext{for the choice of $\sigma$, notice that
		if $|\zeta_1|=|\zeta_2|$,}{}
\begin{equation}
\label{r4}
		-2^{-k}\log \left( \frac{|f_2|}{|f_1|} \right) 
		+2^{-k}\log \left( \frac{|f_1|}{|f_0|} \right) < 2^{-k}\log 4 
\end{equation}
	\par	Thus we may choose $k$ such that $2^{-k} \log 4 < \sigma$.
                In case equation (\ref{r3}) is not satisfied, it is reasonable
		to assume that the two roots have the same modulus indeed.

\section{Numerical Results} 
\label{Sec6}
\begin{figure} 
\vspace{6cm} 
\scalebox{2}{\includegraphics{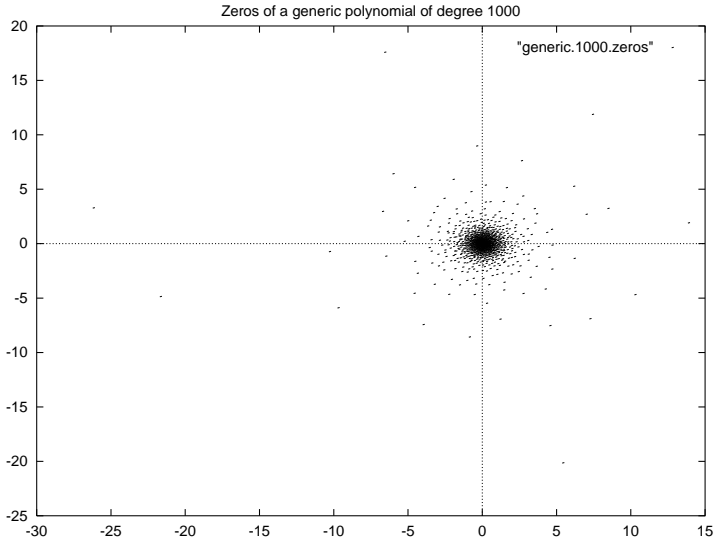}}
\ \\
\caption[ ]{\label{zeros}Zeros of a random polynomial}
\end{figure}
\par
		The results discussed below are tentative, and our
		algorithm deserves further experimentation. Moreover, at
		this moment we are
		using a tentative algorithm to find the argument of the roots.
                Those matters will
		be dealt with in a subsequent paper \newtext{~(\cite{TANGRA})}
		{Written !}. The purpose of this 
		section
		is to illustrate the numerical properties of Renormalized
		Graeffe Iteration, when applied to \oldtext{a generic real
		polynomial} \newtext{random real polynomials}{A3\S 1}.
\par
		A C implementation of our algorithm was tested
                for pseudo-random real and complex
                polynomials. The arguments of the solutions were recovered 
		using an algorithm derived from the Renormalized Graeffe
		Iteration. 
\newtext{	
\par
		IMPLEMENTATION DETAILS: Polynomials with zero coefficients
		do arise in practice, and we explain now how to deal with
		the `log of zero' problem.
		The IEEE floating point arithmetic, implemented by most 
		modern computers, has a few useful features for this sort
		of situation. When one asks a computer to
		produce $\log 0$, it returns a special IEEE value
		called $- \infty$. This is not an error, and further
		calculations can be carried out (as long as they are
		defined). If they are not defined (e.g. $\infty - \infty$)
		then another special value, called a `not a number' is
		returned. 
\par		This last case ($\infty - \infty$) can 
		appear during the computation of a renormalized sum. It
		can be dealt by testing $a$ and $b$ in Example~\ref{ex2}
		for finitesness. In case $a$ or $b$ is infinite, then
		$c$ should be set to $\infty$.
		For an introduction of IEEE arithmetic, see
		~\cite{HIGHAM} pages~45-48 or~\cite{DEMMEL} pages~9-15.
                }{B\S 2(2)}
\oldtext{	
                Results obtained were validated through
		estimates as in~\cite{TCS}. 
		}
		\newtext{The correctness of the results was certified by
		estimates as in~\cite{TCS}. }{A4\S 3(1)}
		\newtext{
		\par
		Experiments were performed in a Pentium 66 
		system running Linux operating system.
		Since the objective here was to illustrate the 
		asymptotic behaviour of the algorithm, we did not
		perform experiments in other systems. Those would be
		necessary if one wanted to compare with other
		algorithms with same asymptotic properties. However,
		this goes far beyond the scope of this paper.}{A4\S 3(2-3)}
\par
		The table below shows the average and median user time
		in a Pentium- based computer, using `long double'
		precision. Time does not include validation time. Ten
		pseudo-random polynomials were tested for each degree. 
		The actual experimental data is plotted in Figure~\ref{Sec6}.
\medskip
\par
\centerline{
\begin{tabular}{||r||r|r||r|r||}
  \hline 
  \hline 
  \multicolumn{1}{||c||}{Degree}&
  \multicolumn{2}{c||}{Real polynomials}&
  \multicolumn{2}{c||}{Complex polynomials} \\
  \hline 
  \multicolumn{1}{||c||}{} & 
  \multicolumn{1}{c}{avg time (s)}&
  \multicolumn{1}{|c||}{median time(s)}&
  \multicolumn{1}{|c}{avg time (s)}&
  \multicolumn{1}{|c||}{median time (s)}\\
  \hline 
  \hline  100 &  0.87 &  0.87 &  1.19 &  1.18 \\
  \hline  200 &  3.23 &  3.24 &  4.35 &  4.30 \\
  \hline  300 &  7.07 &  7.07 &  9.99 &  9.73 \\
  \hline  400 & 12.60 & 12.58 & 17.28 & 17.10 \\
  \hline  500 & 19.41 & 19.38 & 26.75 & 26.71 \\
  \hline  600 & 27.67 & 27.73 & 37.02 & 35.96 \\
  \hline  700 & 37.50 & 37.32 & 51.30 & 49.91 \\
  \hline  800 & 48.89 & 48.72 & 65.39 & 63.61 \\
  \hline  900 & 61.56 & 61.06 & 79.28 & 78.74 \\
  \hline 1000 & 75.89 & 75.47 &102.33 &101.80 \\
  \hline 
  \hline 
\end{tabular}
}
\medskip
\par
	We also computed (approximately) the relative separation of 
        the moduli of the solutions $\zeta_i$:
\[
	\min_{|\zeta_i| < |\zeta_j|}
              \frac{ |\zeta_i| - |\zeta_j| }{\sqrt{1 + |\zeta_i|^2}} 
\]
\par
\begin{figure} 
\vspace{12cm} 
\includegraphics{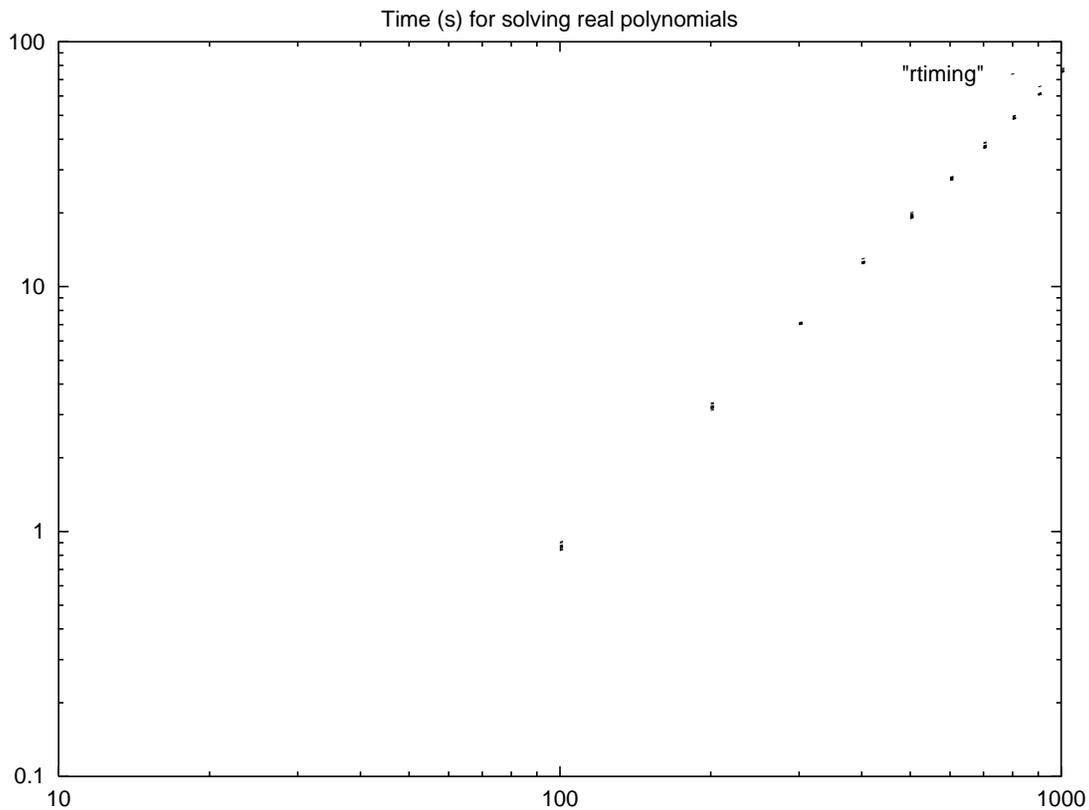}
\ \\
\caption[]{\label{results1}Timing for 100 random real polynomials}
\end{figure}
\begin{figure} 
\vspace{12cm} 
\includegraphics{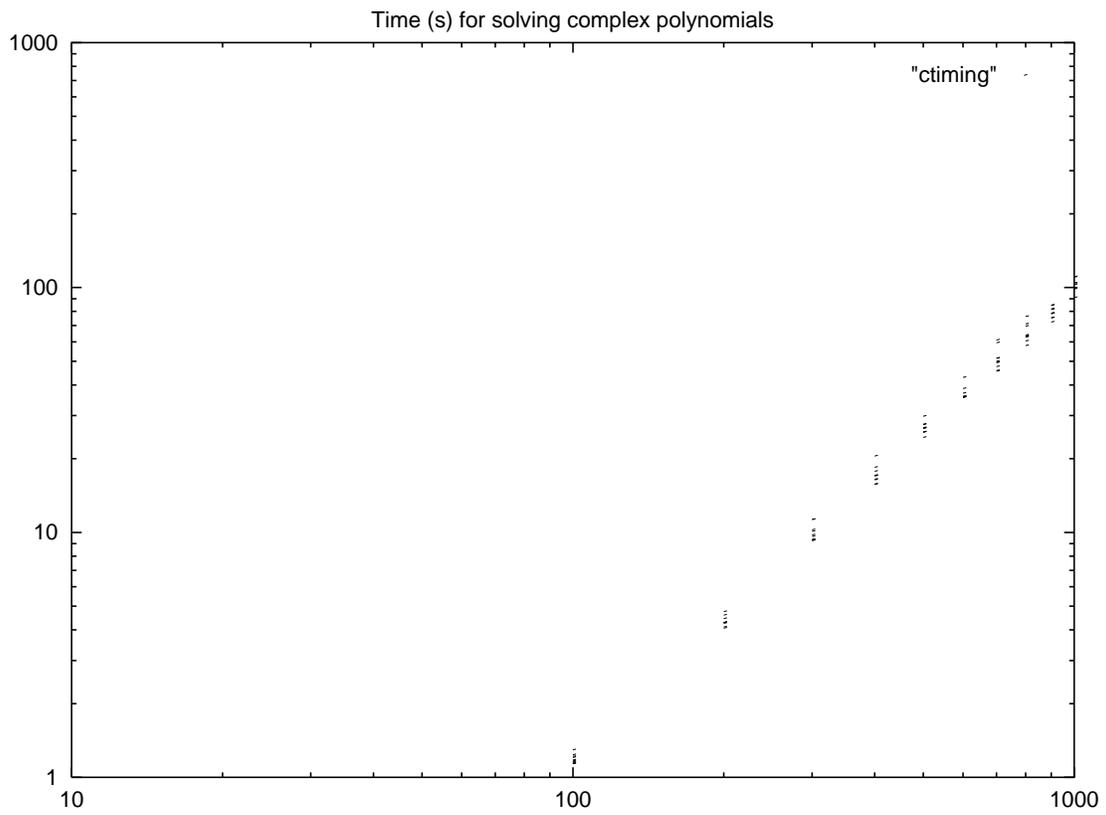}
\ \\
\caption[]{\label{results2}Timing for 100 random complex polynomials}
\end{figure}
\begin{figure} 
\vspace{12cm} 
\includegraphics{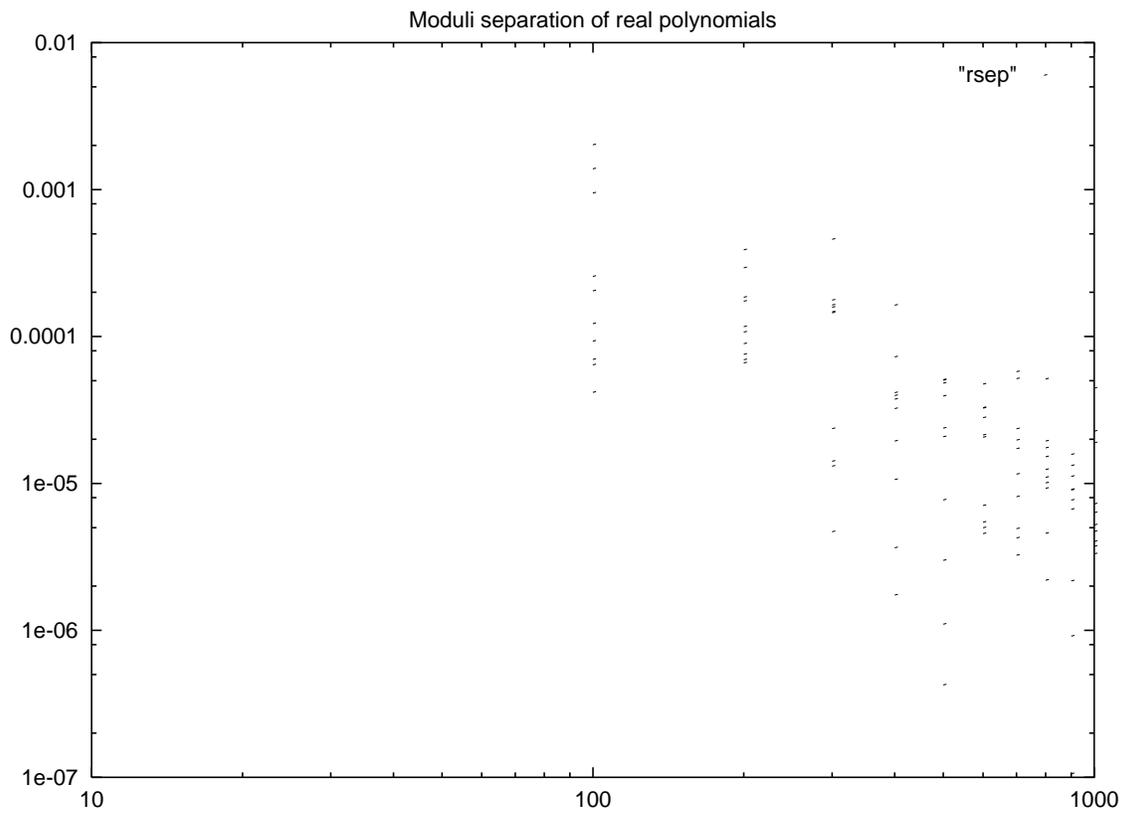}
\ \\
\caption[]{\label{results3}Separation for 100 random real polynomials}
\end{figure}
\begin{figure} 
\vspace{12cm} 
\includegraphics{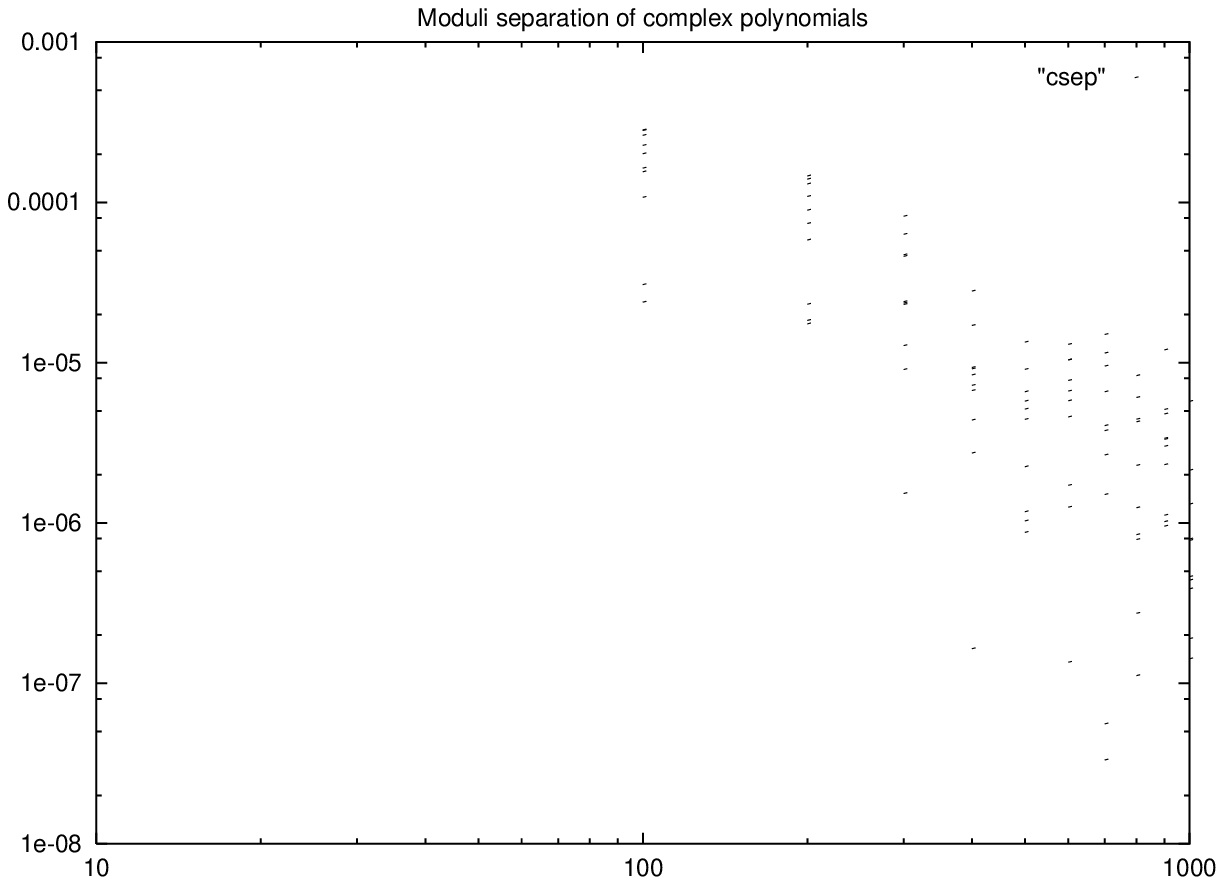}
\ \\
\caption[]{\label{results4}Separation for 100 random complex polynomials}
\end{figure}

	The values obtained are also plotted in 
%	\oldtext{figure~\ref{results}}
	\newtext{figures~\ref{results1} to ~\ref{results4}}{A4\S 3(4)}.	

%\begin{figure} 
%\vspace{8cm} 
%\special{psfile=res1.eps voffset=120 hoffset=-40 vscale=60 hscale=60 }
%\special{psfile=res2.eps voffset=120 hoffset=170 vscale=60 hscale=60 } 
%\special{psfile=res3.eps voffset=-50 hoffset=-40 vscale=60 hscale=60 }
%\special{psfile=res4.eps voffset=-50 hoffset=170 vscale=60 hscale=60 }
%\ \\
%\caption[ ]{\label{results}Time and separation of 100 random polynomials}
% \end{figure}

	\par	Further experimentation is necessary to obtain data
		about polynomials of degree $\gg 1000$. Indeed, due to
		underflow, we cannot represent \oldtext{generic}
		\newtext{random}{A3\S 1}
		high-degree polynomials in the usual floating point
		representation. 
		
\section{Acknowledgments}

JPZ was supported by CNPq and GM was partially supported by CNPq.
The authors wish to thank the following people for their help: 
% suggestions and comments:   
Jean-Pierre Dedieu, Carlos Isnard, Peter Kirrinnis, Welington de Melo, 
Victor Pan, Mike Shub, Jose Felipe da Silva, Steve Smale, Benar Fux Svaiter,
Jean-Claude Yakoubsohn.

The authors gratefully acknowledge a large number of helpful
comments by the anonymous referees, which helped to enhance the
presentation tremendously.

\end{document}